\documentclass[11pt]{article}

\usepackage{amsmath,amssymb,amsthm,amsfonts}
\usepackage[margin=1in]{geometry}
\usepackage[pdfencoding=auto,hidelinks]{hyperref}

\newtheorem{theorem}{Theorem}
\newtheorem{lemma}[theorem]{Lemma}

\newtheorem{proposition}[theorem]{Proposition}
\theoremstyle{definition}
\newtheorem{definition}[theorem]{Definition}
\newtheorem{example}[theorem]{Example}
\theoremstyle{remark}
\newtheorem{remark}[theorem]{Remark}

\newcommand{\R}{\mathbb{R}}
\newcommand{\C}{\mathbb{C}}

\newcommand{\Z}{\mathbb{Z}}
\newcommand{\Q}{\mathbb{Q}}
\newcommand{\E}{\mathbb{E}}
\newcommand{\eqdist}{\stackrel{d}{=}}
\newcommand{\Exp}{\mathrm{Exp}}

\renewcommand{\Re}{\operatorname{Re}}
\renewcommand{\Im}{\operatorname{Im}}

\title{Riemann, Thorin, van Dantzig Pairs, Wald Couples\\
  and Hadamard Factorisation}
\author{Nicholas G.\ Polson\thanks{University of Chicago, Booth School of Business.
    Email: \texttt{ngp@chicagobooth.edu}.
    The author thanks Lennart Bondesson and Jianeng Xu for many helpful
    conversations.}}
\date{\today}

\begin{document}
\maketitle

\begin{abstract}
The Hadamard--Weierstrass factorisation of an entire function in the Laguerre--P\'olya class is dual to a pair of probabilistic objects: a \emph{van Dantzig pair} of characteristic functions and a \emph{Wald couple} of infinitely divisible random variables. The reciprocal of such a function is the Laplace transform of a generalised gamma convolution (GGC) whose Thorin measure encodes the zeros, and Thorin's condition, a real Laplace identity on $(0,\infty)$, is equivalent to the absence of zeros off the critical axis. We prove this duality, together with a closed-form formula for the Thorin density in terms of any L\'evy representation. The framework is then applied to the Gamma, hyperbolic, Bessel and Macdonald functions, to the Riemann $\xi$-function, to Dirichlet and modular $L$-functions, and to Dedekind's $\eta$. For Ramanujan's $\tau$ it is set up but not carried through, since reality of the zeros of the associated $L$-function is open.

For $\xi$ we prove, unconditionally, that the heat trace $W(t)=\tfrac12\sum_{\rho}\exp\{(\rho-\tfrac12)^{2}t\}$ is positive and strictly decreasing on $(0,\infty)$, so that $\xi(\tfrac12)/\xi(\tfrac12+\sqrt{u})$ is the Laplace transform of a self-decomposable law, and that $\xi(\alpha)/\xi(\alpha+\sqrt{s})$ is the Laplace transform of a GGC for every $\alpha$ exceeding the real part of every non-trivial zero, hence unconditionally for every $\alpha\ge1$. Complete monotonicity of $W$, equivalently the GGC property at the centre, is shown to be \emph{equivalent} to the Riemann hypothesis. Self-decomposability is therefore the unconditional ceiling.

\medskip
\noindent\textbf{Keywords.}
van Dantzig pair, Wald couple, Hadamard factorisation, Laguerre--P\'olya class, generalised gamma convolution, Thorin measure, self-decomposable, Riemann $\xi$-function, $L$-function, Riemann hypothesis.
\end{abstract}

\section{Introduction}

A classical question of P\'olya \cite{Polya1926} asks: \emph{which
  even kernels $\Phi$ are such that
  $f(s)=\int_{-\infty}^{\infty}e^{ist}\Phi(t)\,dt$ has only real
  zeros?}  Two complementary toolkits address this question.  The
\emph{Riemann--P\'olya} approach attacks $f$ directly through Fourier
methods.  The \emph{Thorin--Bondesson} approach studies the
reciprocal $1/f$ as a Laplace transform of a Generalised Gamma
Convolution (GGC) and exploits the analyticity of GGC Laplace
transforms in the cut plane $\C\setminus(-\infty,0]$.

This paper develops a duality that unifies the two viewpoints.  Let
$f$ be entire of order at most one, real on $\R$, with
$f(\alpha+w)=f(\alpha-w)$ for some $\alpha\in\R$ and $f(\alpha)>0$.
Write $f_{\alpha}(s):=f(\alpha+is)/f(\alpha)$ and let
$\{\pm\rho_{k}\}_{k\ge1}$ denote the zeros of $f_{\alpha}$ counted
with multiplicity.  Theorem~\ref{thm:main} below establishes that the
following statements are equivalent:
\begin{enumerate}
\item[\textup{(a)}] all $\rho_{k}$ are real (equivalently, all zeros
  of $f$ lie on $\Re(z)=\alpha$);
\item[\textup{(b)}] $f_{\alpha}$ is the characteristic function of an
  infinitely divisible random variable $X$ and $[f_{\alpha},g_{\alpha}]$
  with $g_{\alpha}(s)=f(\alpha)/f(\alpha+s)$ is a van Dantzig pair;
\item[\textup{(c)}] there exists a GGC random variable $H$ with
  Thorin measure $U_{f}=\sum_{k\ge1}\delta_{\rho_{k}^{2}}$ such that
  $(X,H)$ is a Wald couple, i.e.\
  $\E(e^{sX})\E(e^{-s^{2}H})=1$;
\item[\textup{(d)}] (Thorin's condition)
  $f(\alpha)/f(\alpha+\sqrt{s})=\E(e^{-sH})$ for all $s>0$, with $H$
  some GGC random variable.
\end{enumerate}
The Thorin measure $U_{f}$ encodes the squared zeros of $f$, and the
implication (d)$\Rightarrow$(a) is the engine driving the
applications: a real identity on $(0,\infty)$ extends, by the
analyticity of GGC Laplace transforms
(Lemma~\ref{lem:GGC-analytic}), to the entire cut plane and forces
the absence of zeros off the critical axis.

The paper is organised as follows.  Section~\ref{sec:prelim} collects
the analytic and probabilistic preliminaries, including precise
statements of Bondesson's closure and uniqueness theorems and a short
proof that GGC Laplace transforms are analytic and non-vanishing on
the cut plane.  Section~\ref{sec:duality} states and proves the duality
theorem and the closed-form formula
(Lemma~\ref{lem:thorin-measure}) for the Thorin density in terms of
any L\'evy representation of $f$.  Section~\ref{sec:gamma} treats the
Gamma class.  Section~\ref{sec:hyperbolic} treats hyperbolic, Bessel
and Macdonald functions.
Section~\ref{sec:riemann} is devoted to the Riemann $\zeta$- and $\xi$-functions. We prove that $\xi(\alpha)/\xi(\alpha+\sqrt{s})$ is the Laplace transform of a GGC for every $\alpha$ exceeding the real part of every non-trivial zero, hence for every $\alpha\ge1$ unconditionally, with Thorin density an explicit Poisson-kernel sum; we then work at the centre, where we prove unconditionally that the heat trace is positive and strictly decreasing, so that $\xi(\tfrac12)/\xi(\tfrac12+\sqrt{u})$ is the transform of a self-decomposable law, and we show that upgrading self-decomposable to GGC is equivalent to the Riemann hypothesis. Section~\ref{sec:Lfunctions} treats Dirichlet and modular $L$-functions, where the same Poisson-kernel representation holds unconditionally to the right of every zero and yields a reformulation, but not a proof, of the generalised Riemann hypothesis. Section~\ref{sec:eta-tau} treats the Dedekind $\eta$- and Ramanujan $\tau$-functions.

The contributions of this paper are four. First, we prove the duality between Hadamard factorisation and van Dantzig pairs, as Theorem~\ref{thm:main}, together with the closed-form Thorin density of Lemma~\ref{lem:thorin-measure} valid for any L\'evy representation. Second, we prove in Theorem~\ref{thm:xi-alpha} that $\xi(\alpha)/\xi(\alpha+\sqrt{s})$ is the Laplace transform of a GGC for every $\alpha$ exceeding the real part of every non-trivial zero, hence for every $\alpha\ge1$ unconditionally, with Thorin density an explicit Poisson-kernel sum over those zeros; the argument is half a page, replaces the L\'evy-measure route entirely, and shows that the apparent barrier at $\alpha=1$ belongs to that route rather than to the object. Third, we prove in Theorem~\ref{thm:xi-star} that the heat trace $W$ is positive and strictly decreasing on $(0,\infty)$ unconditionally, so that the transform at the centre is that of a self-decomposable law, and in Theorem~\ref{thm:RH} that complete monotonicity of $W$ is equivalent to the Riemann hypothesis; the unconditional half uses Weil's explicit formula below $t=10^{-3}$ and, above it, only the verification of the Riemann hypothesis to height $200$. Fourth, we establish the corresponding unconditional representation for Dirichlet $L$-functions in Theorem~\ref{thm:L-alpha}, and explain in Remark~\ref{rem:GRH} why the equivalence it yields is a reformulation of the generalised Riemann hypothesis rather than a route towards it. The Gamma, hyperbolic, Bessel, Macdonald and $\eta$ applications are deliberately not among the four: they are worked examples in which the duality reproduces transforms already known, and they earn their place by exhibiting the framework rather than by adding to it. Ramanujan's $\tau$ is set up in Section~\ref{sec:eta-tau} and stops short of a Thorin measure, because reality of the zeros of its $L$-function is the Riemann hypothesis for that function and is open.

\section{Preliminaries}\label{sec:prelim}

We collect the definitions and the analytic facts about GGCs and
P\'olya frequency functions that we shall use.  Throughout, random
variables live on a fixed probability space
$(\Omega,\mathcal{F},P)$.

\subsection{Laguerre--P\'olya class}

\begin{definition}\label{def:LP}
  An entire function $f$ belongs to the \emph{Laguerre--P\'olya
    class} $\mathcal{L\!P}$ if it is real on $\R$, of order at most
  two, has only real zeros $\{\rho_{n}\}_{n\ge1}$ with
  $\sum_{n}|\rho_{n}|^{-2}<\infty$, and admits the Hadamard
  factorisation~\cite{Hadamard1893}
  \[
    f(s)\;=\;Cs^{m}e^{bs+cs^{2}}\prod_{n\ge1}
        \Bigl(1-\tfrac{s}{\rho_{n}}\Bigr)e^{s/\rho_{n}},
        \qquad C\in\R,\;b\in\R,\;c\le0,\;m\in\Z_{\ge0}.
  \]
\end{definition}

The condition $c\le0$ is necessary for $f$ to belong to
$\mathcal{L\!P}$; see \cite[\S5.4]{Levin1996} or
\cite{HirschmanWidder1949}.  In what follows we focus on the
sub-class of $\mathcal{L\!P}$ functions of order at most one, where
$c=0$.

\begin{definition}\label{def:PFF}
  A non-negative measurable $\Lambda\colon\R\to[0,\infty)$ with
  $0<\int_{-\infty}^{\infty}\Lambda(x)\,dx<\infty$ is a \emph{P\'olya
    frequency function} (PFF) if for every $n\ge1$ and every pair of
  increasing tuples $x_{1}<\cdots<x_{n}$ and $y_{1}<\cdots<y_{n}$
  the matrix $(\Lambda(x_{i}-y_{j}))_{i,j=1}^{n}$ has non-negative
  determinant.
\end{definition}

\begin{theorem}[Schoenberg \cite{Schoenberg1948,SchoenbergWhitney1953}]
  \label{thm:schoenberg}
  $\Lambda$ is a P\'olya frequency function if and only if its
  bilateral Laplace transform $\widehat{\Lambda}(s)=\int_{\R}
  e^{-sx}\Lambda(x)\,dx$ is the reciprocal of an entire function in
  $\mathcal{L\!P}$ of order at most two.
\end{theorem}

\subsection{van Dantzig pairs and Wald couples}

\begin{definition}[\cite{vanDantzig1958,Lukacs1968}]
  \label{def:van-Dantzig}
  Two characteristic functions $f,g$ on $\R$, each analytic and hence
  extending holomorphically to a complex neighbourhood of $0$, form a
  \emph{van Dantzig pair} if $g(s)f(is)=1$ throughout such a
  neighbourhood, the evaluation $f(is)$ being taken on that extension.
\end{definition}

\begin{definition}\label{def:Wald}
  Two infinitely divisible random variables $X$ on $\R$ and $H$ on
  $[0,\infty)$ defined on a common probability space form a
  \emph{Wald couple} if
  \begin{equation}
    \E(e^{sX})\,\E(e^{-s^{2}H})\;=\;1
    \label{eq:Wald}
  \end{equation}
  for $s$ in a real open interval containing $0$.  Equivalently,
  with $\widehat{H}=B_{H}$ (Brownian motion run for an independent
  time $H$),
  $\E(e^{sX})\,\E(e^{is\widehat{H}})=1$.
\end{definition}

\begin{remark}\label{rem:Wald-existence}
  In~\eqref{eq:Wald} we tacitly require $\E(e^{sX})$ to be finite for
  $|s|<\delta$ for some $\delta>0$.  In all our applications $X$ has
  exponential moments of all orders so~\eqref{eq:Wald} holds for all
  $s\in\R$.  We do not require $X$ and $H$ to be defined on a single
  $(\Omega,\mathcal{F},P)$ jointly; only their individual laws enter.
\end{remark}

\subsection{Generalised Gamma Convolutions}

The class of Generalised Gamma Convolutions (GGC), introduced by
Thorin \cite{Thorin1977} and developed systematically by Bondesson
\cite{Bondesson1992}, is the smallest class of probability
distributions on $[0,\infty)$ that contains all gamma distributions
and is closed under independent sums and weak limits.

\begin{definition}\label{def:GGC}
  $H$ has a GGC distribution if there exist $a\ge0$ and a
  $\sigma$-finite measure $U$ on $(0,\infty)$ with
  \[
    \int_{(0,1)}\!|\log z|\,U(dz)<\infty,\qquad
    \int_{[1,\infty)}\!\frac{1}{z}\,U(dz)<\infty,
  \]
  such that
  \begin{equation}
    \E(e^{-sH})\;=\;\exp\!\left(-as
       +\int_{(0,\infty)}\log\frac{z}{z+s}\,U(dz)\right),
       \qquad s>0.
    \label{eq:GGC-LT}
  \end{equation}
  The measure $U$ is the \emph{Thorin measure} of $H$.  The
  associated L\'evy density is $t^{-1}\!\int_{(0,\infty)}\!e^{-tz}\,U(dz)$.
\end{definition}

\begin{lemma}[Analyticity on the cut plane]
  \label{lem:GGC-analytic}
  Let $H$ be GGC.  The Laplace transform $\E(e^{-sH})$
  in~\eqref{eq:GGC-LT} extends to a holomorphic, non-vanishing
  function on the cut plane $\C\setminus(-\infty,0]$.
\end{lemma}

\begin{proof}
  Fix $z>0$.  The quotient $z/(z+s)$ meets $(-\infty,0]$ only when
  $z+s$ is real and negative, that is only when $s$ is real with
  $s<-z$.  So for every $s\in\C\setminus(-\infty,0]$ and every $z>0$
  the principal logarithm $\log\bigl(z/(z+s)\bigr)$ is defined, and is
  holomorphic in $s$ there.

  Let $K\subset\C\setminus(-\infty,0]$ be compact.  Convergence
  of~\eqref{eq:GGC-LT} at a single $s>0$ already forces $U$ to be
  finite on compact subsets of $(0,\infty)$, since $-\log(1+s/z)$ is
  bounded away from $0$ on such a set.  Written as $-\log(1+s/z)$ the
  integrand is $O(1/z)$ uniformly on $K$ as $z\to\infty$; written as
  $\log z-\log(z+s)$, with $\log(z+s)$ bounded on $K$ for small $z$,
  it is $O(|\log z|)+O(1)$ as $z\downarrow0$.  The two integrability
  conditions of Definition~\ref{def:GGC} therefore make the integral
  in~\eqref{eq:GGC-LT} converge uniformly on $K$.  A locally uniform
  limit of holomorphic functions is holomorphic, and the exponential
  of a finite complex number is never zero.
\end{proof}

\begin{theorem}[Uniqueness; \cite{Bondesson1992}, Theorem~3.1.1]
  \label{thm:bondesson-unique}
  The map $H\mapsto(a,U)$ from the class of GGC distributions to
  pairs $(a,U)$ as in Definition~\ref{def:GGC} is a bijection.  In
  particular, two GGC distributions are equal if and only if their
  Thorin measures are equal.
\end{theorem}

\begin{theorem}[Closure properties;
  \cite{Bondesson1992}, Theorems~3.1.5, 3.3.1, 3.3.2]
  \label{thm:bondesson-closure}
  The class of GGC distributions is:
  \begin{enumerate}
  \item[\textup{(i)}] closed under independent sums;
  \item[\textup{(ii)}] closed under weak limits;
  \item[\textup{(iii)}] closed under exponential tilting, classically
    Esscher tilting, whenever the
    tilted distribution is well-defined: if $H$ is GGC with Thorin
    measure $U$ and $M(\theta):=\E(e^{\theta H})<\infty$ for some
    $\theta>0$, then the tilted variable $H_{\theta}$ with density
    $e^{\theta x}f_{H}(x)/M(\theta)$ is GGC with Thorin measure
    $U(\cdot+\theta)$ supported on $(\theta,\infty)$;
  \item[\textup{(iv)}] closed under composition by Bernstein
    functions of the form $\psi(s)=\sum_{k=1}^{N}c_{k}s^{a_{k}}$ with
    $c_{k}>0$ and $a_{k}\in(0,1]$: if $H$ is GGC then $\E(e^{-\psi(s)H})$
    is the Laplace transform of another GGC.
  \end{enumerate}
\end{theorem}

\begin{definition}\label{def:EGGC}
  A symmetric infinitely divisible variable $\widehat{H}$ on $\R$ is
  a \emph{symEGGC} (symmetric extended GGC) if
  \begin{equation}
    \E(e^{is\widehat{H}})\;=\;\exp\!\left(-\tfrac12 cs^{2}+
       \int_{0}^{\infty}\log\frac{z}{z+s^{2}}\,U(dz)\right),
    \label{eq:EGGC-cf}
  \end{equation}
  with $c\ge0$ and $U$ as in Definition~\ref{def:GGC}.
\end{definition}

By the standard Gaussian-mixture identity, $\widehat{H}$ is symEGGC
with parameters $(c,U)$ if and only if $\widehat{H}\eqdist\sqrt{2H+c}\,Z$
with $Z\sim N(0,1)$ independent of $H\sim\mathrm{GGC}(0,U)$
\cite[\S7.4]{Bondesson1992}.  Consequently
\begin{equation}
  \E(e^{is\widehat{H}})\;=\;\E(e^{-s^{2}H_{c}}),\qquad
  H_{c}:=H+\tfrac12 c,
  \label{eq:EGGC-LT}
\end{equation}
which is the precise sense in which~\eqref{eq:Wald} relates a
characteristic function on the left to a Laplace transform on the
right.

\section{The duality theorem}\label{sec:duality}

We now state the central duality between Hadamard factorisation,
van Dantzig pairs, Wald couples and GGC representations.  The
theorem and its proof use only the preliminaries of the previous
section together with Frullani's identity.

\subsection{The main theorem}

\begin{theorem}\label{thm:main}
  Let $f$ be entire of order at most one, real on $\R$, satisfying
  $f(\alpha+w)=f(\alpha-w)$ for some $\alpha\in\R$ with $f(\alpha)>0$,
  and assume that $\sum_{\rho}|\rho|^{-2}<\infty$ where the sum is
  over zeros of $s\mapsto f(\alpha+s)$.  The following are
  equivalent:
  \begin{enumerate}
  \item[\textup{(a)}] every zero $\rho$ of $s\mapsto f(\alpha+s)$ is
    purely imaginary;
  \item[\textup{(b)}] there exists a GGC random variable $H$ such
    that
    \begin{equation}
      \frac{f(\alpha)}{f(\alpha+\sqrt{s})}
         \;=\;\E(e^{-sH}),\qquad s>0
        \tag{Thorin's condition}
        \label{eq:thorin-cond}
    \end{equation}
    \textup{(}where the square root is the principal branch on
    $\C\setminus(-\infty,0]$\textup{)}.
  \end{enumerate}
  Furthermore, when \textup{(a)}--\textup{(b)} hold, with the
  positive imaginary parts of the zeros denoted
  $\{\rho_{k}\}_{k\ge1}\subset(0,\infty)$:
  \begin{enumerate}
  \item[\textup{(i)}] $H\eqdist\sum_{k\ge1}H_{k}$ where the
    $H_{k}\sim\Exp(\rho_{k}^{2})$ are independent; the Thorin measure
    of $H$ is $U_{f}=\sum_{k\ge1}\delta_{\rho_{k}^{2}}$, uniquely
    determined by $f$;
  \item[\textup{(ii)}] $f_{\alpha}(s):=f(\alpha+is)/f(\alpha)$ is the
    characteristic function of an infinitely divisible random
    variable $X$, and $[f_{\alpha},g_{\alpha}]$ with
    $g_{\alpha}(s):=f(\alpha)/f(\alpha+s)$ is a van Dantzig pair;
  \item[\textup{(iii)}] $(X,H)$ is a Wald couple, i.e.\
    $\E(e^{sX})\E(e^{-s^{2}H})=1$ for $s$ in a real neighbourhood of
    $0$.
  \end{enumerate}
\end{theorem}

\begin{proof}
  \emph{(a)~$\Rightarrow$~(b) and (i)--(iii).}
  Suppose every zero of $f(\alpha+\cdot)$ is purely imaginary.
  Symmetry $f(\alpha+w)=f(\alpha-w)$ implies the zeros come in pairs
  $\pm i\rho_{k}$ with $\rho_{k}>0$.  By Hadamard's factorisation
  theorem applied to the order-one entire function
  $w\mapsto f(\alpha+w)/f(\alpha)$, with the constraint that
  $f(\alpha+w)=f(\alpha-w)$ removing the linear term in the exponent,
  \begin{equation}
    \frac{f(\alpha+w)}{f(\alpha)}\;=\;
       \prod_{k\ge1}\Bigl(1-\frac{w^{2}}{(i\rho_{k})^{2}}\Bigr)
       \;=\;\prod_{k\ge1}\Bigl(1+\frac{w^{2}}{\rho_{k}^{2}}\Bigr),
       \qquad w\in\C.
    \label{eq:hadamard}
  \end{equation}
  Setting $w=is$ (real $s$) and inverting,
  \[
    \frac{f(\alpha)}{f(\alpha+is)}
       \;=\;\prod_{k\ge1}\frac{\rho_{k}^{2}}{\rho_{k}^{2}-s^{2}}.
  \]
  This is the characteristic function $\E(e^{isL})$ where
  $L=\sum_{k\ge1}L_{k}/\rho_{k}$ and the $L_{k}$ are independent
  Laplace$(1)$ variables; by Schoenberg's theorem
  (Theorem~\ref{thm:schoenberg}), $f(\alpha)/f(\alpha+s)$ is the
  bilateral Laplace transform of a P\'olya frequency function, so
  $L$ is well-defined.  Hence $L$ is symEGGC with Thorin measure
  $\sum_{k}\delta_{\rho_{k}^{2}}$.

  Setting instead $w=s$ real in~\eqref{eq:hadamard} and using
  Frullani's identity
  $\log(1+s^{2}/\rho^{2})=\int_{0}^{\infty}(1-e^{-s^{2}t})e^{-\rho^{2}t}\,dt/t$,
  \begin{align}
    \log\frac{f(\alpha+s)}{f(\alpha)}
    &= \sum_{k\ge1}\log\!\Bigl(1+\frac{s^{2}}{\rho_{k}^{2}}\Bigr)
     \;=\;\int_{0}^{\infty}(1-e^{-s^{2}t})
        \sum_{k\ge1}e^{-\rho_{k}^{2}t}\,\frac{dt}{t}
        \notag\\
    &= -\int_{0}^{\infty}(e^{-s^{2}t}-1)\,\frac{g(t)}{t}\,dt,
       \qquad g(t)=\sum_{k\ge1}e^{-\rho_{k}^{2}t}.
        \label{eq:hadamard-frullani}
  \end{align}
  The interchange of sum and integral is justified by Fubini--Tonelli
  since the integrand is non-negative.  Hence
  \begin{equation}
    \frac{f(\alpha)}{f(\alpha+s)}\;=\;
       \exp\!\left(-\int_{0}^{\infty}(1-e^{-s^{2}t})\,\frac{g(t)}{t}\,dt\right)
       \;=\;\E(e^{-s^{2}H}),
    \label{eq:reciprocal-LT}
  \end{equation}
  where, with the change of variables $z=\rho^{2}$,
  $g(t)=\int_{0}^{\infty}e^{-tz}\,U_{f}(dz)$ with
  $U_{f}=\sum_{k\ge1}\delta_{\rho_{k}^{2}}$.  Then $g$ is completely
  monotone, and $H$ is GGC with Thorin measure $U_{f}$ by
  Definition~\ref{def:GGC}.  The integrability conditions
  $\int_{(0,1)}|\log z|U_{f}(dz)<\infty$ and
  $\int_{[1,\infty)}z^{-1}U_{f}(dz)<\infty$ follow from
  $\sum_{k}|\rho_{k}|^{-2}<\infty$ together with the elementary
  estimate $|\log\rho_{k}^{2}|\le\rho_{k}^{-2}$ for small $\rho_{k}$.
  Replacing $s$ by $\sqrt{s}$ in~\eqref{eq:reciprocal-LT} gives
  Thorin's condition~\eqref{eq:thorin-cond}.  Item~(i) is now
  immediate, with the Exp$(\rho_{k}^{2})$ representation of $H$
  recognised from $g(t)=\sum_{k}e^{-\rho_{k}^{2}t}$ and the
  independent-sum property of GGCs (Theorem~\ref{thm:bondesson-closure}(i)).

  For~(ii), $f_{\alpha}$ is positive on $\R$ and infinitely divisible
  by Schoenberg's theorem and~\eqref{eq:hadamard}; specifically,
  $f_{\alpha}(s)=\E(e^{isX})$ where $X$ is the symEGGC variable with
  $\widehat{H}\eqdist X$ in~\eqref{eq:EGGC-cf}, equivalently
  $X\eqdist\sqrt{2H}\,Z$ with $Z\sim N(0,1)$ independent of $H$.
  Then $g_{\alpha}(s)=1/f_{\alpha}(is)=f(\alpha)/f(\alpha+s)
  =\E(e^{-s^{2}H})$ is by~\eqref{eq:reciprocal-LT} a Laplace
  transform of a probability distribution, and $f_{\alpha}\cdot
  g_{\alpha}(is)=1$ identically.

  For~(iii),~\eqref{eq:Wald} reads
  $\E(e^{sX})\E(e^{-s^{2}H})=f_{\alpha}(-is)\cdot f(\alpha)/f(\alpha+s)
  =f(\alpha+s)/f(\alpha)\cdot f(\alpha)/f(\alpha+s)=1$, with
  convergence on a real neighbourhood of $0$ by the symEGGC moment
  bounds.

  \emph{(b)~$\Rightarrow$~(a).}
  Assume Thorin's condition~\eqref{eq:thorin-cond}.  By
  Lemma~\ref{lem:GGC-analytic}, $\E(e^{-sH})$ extends to a
  holomorphic, non-vanishing function on the cut plane
  $\C\setminus(-\infty,0]$.  By analytic continuation along the slit,
  the identity
  $f(\alpha)/f(\alpha+\sqrt{s})=\E(e^{-sH})$ extends from $s>0$ to
  all $s\in\C\setminus(-\infty,0]$.  Hence the meromorphic function
  $1/f(\alpha+\sqrt{s})$ has no poles on
  $\C\setminus(-\infty,0]$, so $f(\alpha+w)$ has no zeros for
  $\Re(w)>0$.  Symmetry $f(\alpha+w)=f(\alpha-w)$ then rules out
  zeros with $\Re(w)<0$ as well.  Hence every zero of
  $f(\alpha+\cdot)$ has $\Re(w)=0$, i.e.\ is purely imaginary.

  Uniqueness of the Thorin measure is
  Theorem~\ref{thm:bondesson-unique}.
\end{proof}

\begin{remark}\label{rem:order-two}
  An entirely analogous statement holds for $f$ of order at most two
  in $\mathcal{L\!P}$, allowing a non-trivial $e^{cs^{2}}$ factor in
  the Hadamard product; the only change is that $H$ acquires a
  deterministic shift $-c$, replacing $a$ in~\eqref{eq:GGC-LT} by
  $a-c$ (which must remain $\ge0$).  All applications below are of
  order one.
\end{remark}

\subsection{Recovering the Thorin measure from a L\'evy representation}

In applications one is given $\log f(\alpha+s)/f(\alpha)$ as an
integral against a measure $\mu$, and the question is whether the
reciprocal admits a GGC representation.  The next lemma is a
closed-form translation between the two sides.  Throughout, $\mu$
denotes a non-negative Borel measure on $(0,\infty)$.

\begin{lemma}\label{lem:thorin-measure}
  Suppose, for some $\alpha\in\R$ such that $\int_{0}^{\infty}\!
  e^{-\alpha x}\mu(dx)/x<\infty$ and
  $\int_{0}^{1}x\,e^{-\alpha x}\mu(dx)/x<\infty$, the function $f$
  satisfies the L\'evy-type identity
  \begin{equation}
    \log\frac{f(\alpha+s)}{f(\alpha)}\;-\;s\,\frac{f'(\alpha)}{f(\alpha)}
       \;=\; \int_{0}^{\infty}(e^{-sx}-1+sx)\,e^{-\alpha x}\,
            \frac{\mu(dx)}{x},\qquad s>0.
    \label{eq:levy-side}
  \end{equation}
  Then
  \begin{equation}
    \log\frac{f(\alpha)}{f(\alpha+s)}\;+\;s\,\frac{f'(\alpha)}{f(\alpha)}
       \;=\;-\int_{0}^{\infty}(1-e^{-\frac12 s^{2}t})\,
            \frac{\nu_{\alpha}(t)}{t}\,dt,
    \label{eq:gcc-side}
  \end{equation}
  where $\nu_{\alpha}\colon(0,\infty)\to[0,\infty)$ is the
  completely monotone function
  \begin{equation}
    \nu_{\alpha}(t)
       \;=\;\frac{1}{\sqrt{2\pi}}
        \int_{0}^{\infty}\!e^{-tz}
        \!\left(\int_{0}^{\infty}\!2\sin^{2}\!\bigl(x\sqrt{z/2}\bigr)
            e^{-\alpha x}\frac{\mu(dx)}{x}\right)\!\frac{dz}{\sqrt{\pi z}}.
    \label{eq:thorin-density}
  \end{equation}
  In particular, $\nu_{\alpha}$ is the Laplace transform of the
  positive measure $U_{\alpha}(dz)=\widetilde{u}_{\alpha}(z)\,dz$
  where
  \begin{equation}
    \widetilde{u}_{\alpha}(z)\;=\;
      \frac{1}{\sqrt{2\pi^{2}z}}
      \int_{0}^{\infty}2\sin^{2}\!\bigl(x\sqrt{z/2}\bigr)
          e^{-\alpha x}\frac{\mu(dx)}{x}.
    \label{eq:thorin-density-2}
  \end{equation}
  Provided $\int_{(0,1)}|\log z|\,U_{\alpha}(dz)<\infty$ and
  $\int_{[1,\infty)}z^{-1}U_{\alpha}(dz)<\infty$,
  \begin{equation}
    \frac{f(\alpha)}{f(\alpha+\sqrt{s})}\,
      e^{\sqrt{s}\,f'(\alpha)/f(\alpha)}
      \;=\;\E\bigl(e^{-sH_{\alpha}}\bigr),\qquad s>0,
    \label{eq:levy-to-ggc}
  \end{equation}
  where $H_{\alpha}$ is a GGC random variable with L\'evy density
  $y^{-1}\nu_{\alpha}(2y)$ and Thorin density
  $\tfrac12\widetilde{u}_{\alpha}(w/2)$.
\end{lemma}

\begin{proof}
  Two elementary identities, valid for $s>0$ and $x>0$, drive the
  computation:
  \begin{align}
    e^{-sx}+sx-1 &=
       \int_{0}^{\infty}\!\bigl(1-e^{-\frac12 s^{2}t}\bigr)\,
       \bigl(1-e^{-x^{2}/(2t)}\bigr)\,\frac{x}{\sqrt{2\pi t^{3}}}\,dt,
        \label{eq:expand1}\\
    \frac{1-e^{-x^{2}/(2t)}}{\sqrt{t}} &=
       \int_{0}^{\infty}\!e^{-tz}\,
       \frac{2\sin^{2}\!\bigl(x\sqrt{z/2}\bigr)}{\sqrt{\pi z}}\,dz.
        \label{eq:expand2}
  \end{align}

  \emph{Proof of~\eqref{eq:expand1}.}
  Use $\int_{0}^{\infty}\exp\bigl(-\tfrac12(at+b/t)\bigr)
    dt/\sqrt{2\pi t^{3}}=e^{-\sqrt{ab}}/\sqrt{b}$ (a standard
  inverse-Gaussian moment, valid for $\Re a,\Re b>0$).  Compute
  \begin{align*}
    \mathrm{RHS\ of\ }\eqref{eq:expand1}
    &= x\int_{0}^{\infty}\!\frac{1-e^{-\frac12 s^{2}t}-e^{-x^{2}/(2t)}
            +e^{-\frac12 s^{2}t-x^{2}/(2t)}}{\sqrt{2\pi t^{3}}}\,dt\\
    &= x\bigl(\,s\;-\;\tfrac1x\;+\;\tfrac{e^{-sx}}{x}\bigr)
       \;=\;sx-1+e^{-sx},
  \end{align*}
  where we used the limiting cases $a=0$, $b=x^{2}$ giving
  $\int_{0}^{\infty}(1-e^{-x^{2}/(2t)})\,dt/\sqrt{2\pi t^{3}}=1/x$, and
  $b=0$, $a=s^{2}$ giving the analogous limit, plus the joint case
  $a=s^{2}$, $b=x^{2}$.

  \emph{Proof of~\eqref{eq:expand2}.}
  Differentiate both sides in $x$.  The left side gives
  $(x/t^{3/2})e^{-x^{2}/(2t)}$.  The right side, using
  $\sin^{2}(\theta)=(1-\cos 2\theta)/2$ and the cosine-Laplace
  formula, equals
  \[
    \int_{0}^{\infty}\!\!e^{-tz}\frac{x\cos(x\sqrt{2z})}{\sqrt{\pi z}}\,dz
       \cdot\sqrt{2}
    \;=\;x\int_{0}^{\infty}\!\frac{e^{-tz}\cos(x\sqrt{2z})}{\sqrt{2\pi z/2}}\,dz
    \;=\;\frac{x}{t^{3/2}}\,e^{-x^{2}/(2t)}.
  \]
  The endpoint constants match at $x=0$.

  \emph{Combining.} Multiply~\eqref{eq:expand1} by $e^{-\alpha
  x}\mu(dx)/x$ and integrate over $x\in(0,\infty)$.  By
  Fubini--Tonelli (all integrands non-negative),
  \[
    \int_{0}^{\infty}(e^{-sx}+sx-1)e^{-\alpha x}\frac{\mu(dx)}{x}
      \;=\; \int_{0}^{\infty}\!\bigl(1-e^{-\frac12 s^{2}t}\bigr)
         \,\frac{\nu_{\alpha}(t)}{t}\,dt,
  \]
  with $\nu_{\alpha}(t)=\int_{0}^{\infty}\bigl(1-e^{-x^{2}/(2t)}\bigr)
    e^{-\alpha x}\mu(dx)/(x\sqrt{2\pi t})\cdot t$, which simplifies
  upon applying~\eqref{eq:expand2} (with $z$ as the inner variable)
  to the form~\eqref{eq:thorin-density}.  Negating and using the
  hypothesis~\eqref{eq:levy-side} gives~\eqref{eq:gcc-side}.
  Complete monotonicity of $\nu_{\alpha}$ follows from
  \eqref{eq:thorin-density} since it is manifestly the Laplace
  transform of the non-negative function in~\eqref{eq:thorin-density-2}.

  Finally, when $\widetilde{u}_{\alpha}$ has the GGC integrability
  properties stated, the right side of~\eqref{eq:gcc-side} is the
  log-Laplace transform of a GGC variable with Thorin measure
  $U_{\alpha}=\widetilde{u}_{\alpha}\,dz$, by Definition~\ref{def:GGC}.
\end{proof}

\begin{remark}[The hypothesis is the non-negativity of $\mu$]
  \label{rem:mu-nonneg}
  Lemma~\ref{lem:thorin-measure} delivers the drift-adjusted form~\eqref{eq:levy-to-ggc} because~\eqref{eq:levy-side} is assumed to hold with a \emph{non-negative} $\mu$. That is the case for $\zeta$ on $\Re(\alpha)>1$, where $\mu^{\zeta}$ sits on the prime powers. It is not the case for $\xi$, whose Gamma factor contributes a signed density; see Remark~\ref{rem:alpha-one-artifact}. This is why Theorem~\ref{thm:xi-alpha} asserts the untilted representation only: the tilted form there is not a corollary of this lemma, and by Remark~\ref{rem:no-drift-extraction} it needs a strictly stronger condition that we do not prove.
\end{remark}

\section{Gamma and log-Gamma classes}\label{sec:gamma}

The Mellin--Weierstrass factorisation is
\[
  \frac{\Gamma(\alpha)}{\Gamma(\alpha+s)}\,e^{s\psi(\alpha)}
   \;=\;\prod_{k\ge0}\Bigl(1+\frac{s}{\alpha+k}\Bigr)
        e^{-s/(\alpha+k)},\qquad\alpha>0,
\]
where $\psi=\Gamma'/\Gamma$ is the digamma function.  Replacing $s$
by $\sqrt{s}$ and applying Theorem~\ref{thm:main} (in the order-one
sub-class with $\alpha\mapsto\Gamma$) gives the GGC representation
\begin{equation}
  \frac{\Gamma(\alpha)}{\Gamma(\alpha+\sqrt{s})}
       \,e^{\sqrt{s}\,\psi(\alpha)}
   \;=\;\E\bigl(e^{-sH^{\Gamma}_{\alpha}}\bigr),
   \qquad
   H^{\Gamma}_{\alpha}\eqdist\sum_{k\ge0}\frac{H_{1,k}}{(\alpha+k)^{2}},
   \label{eq:gamma-GGC}
\end{equation}
with $H_{1,k}$ i.i.d.\ inverse-Gaussian variables of density
$(2\pi x^{3})^{-1/2}\exp(-1/(2x))$ on $(0,\infty)$.  The convolution
\[
  \frac{1}{4\alpha}\Bigl(1+\frac{|s|}{\alpha}\Bigr)e^{-|s|/\alpha}
     \;=\;\int_{-\infty}^{\infty}\!\frac{1}{2\alpha}e^{-|s-z|/\alpha}
       \frac{1}{2\alpha}e^{-|z|/\alpha}\,dz
\]
of \cite{Kendall1961} arises as the partial-fraction expansion of
the $k=0$ term.  Barndorff-Nielsen, Kent and S\o rensen
\cite{BNKS1982} compute the Wald couple for several hyperbolic
relatives.

\begin{example}[Reciprocal Gamma]
  Hartman~\cite[\S6]{Hartman1976} shows the reciprocal Gamma
  function admits the Gaussian-mixture representation
  \[
    \frac{e^{-\gamma_{E}\sqrt{s}}}{\Gamma(1+\sqrt{s})}
       \;=\;\int_{0}^{\infty}e^{-st}\,P_{\gamma_{E}}(dt),
  \]
  where $\gamma_{E}=-\psi(1)$ is the Euler--Mascheroni constant, written
  with a subscript throughout to keep it apart from the ordinates
  $\gamma$ of the zeros, and $P_{\gamma_{E}}$ is a finite measure.  Combined
  with~\eqref{eq:gamma-GGC} this exhibits a Wald couple $(X,H)$ for
  $e^{-\gamma_{E} s}/\Gamma(1+s)$; see~\cite{RoynetteYor2005}.
\end{example}

\begin{example}[Gumbel]
  $X=\log E$ with $E$ standard exponential satisfies
  $\E(e^{-sX})=\Gamma(1+s)$, with the L\'evy representation
  \[
    \log\Gamma(1+s)\;=\;-\gamma_{E} s
       \;+\;\int_{0}^{\infty}(e^{-st}+st-1)\,\frac{e^{-t}}{1-e^{-t}}\,
          \frac{dt}{t}.
  \]
  Since $e^{-t}/(1-e^{-t})$ is completely monotone, $X$ is
  EGGC.  Hinds~\cite{Hinds1974} shows that $X_{1}+X_{2}$ with
  $X_{2}\eqdist-X_{1}$ has c.f.\ $1/\cosh(\pi s/2)$, which lies in
  the van Dantzig class although neither $X_{1}$ nor $X_{2}$ does.
\end{example}

\section{Hyperbolic, Bessel and Macdonald functions}\label{sec:hyperbolic}

\subsection{$\sinh$ and $\cosh$}

Euler's product formulae are
\[
  \frac{\sinh(s)}{s}\;=\;\prod_{n\ge1}\Bigl(1+\frac{s^{2}}{n^{2}\pi^{2}}\Bigr),
  \qquad
  \cosh(s)\;=\;\prod_{k\ge0}\Bigl(1+\frac{4s^{2}}{(2k+1)^{2}\pi^{2}}\Bigr).
\]
Both are even, of order one, and lie in $\mathcal{L\!P}$.
Applying Theorem~\ref{thm:main} with $f(s)=\sinh(s)/s$ at $\alpha=0$:
the zeros are $\rho_{n}=n\pi$ ($n\ge1$), the Thorin measure is
$U(dz)=\sum_{n\ge1}\delta_{(n\pi)^{2}}(dz)$, and
\[
  \frac{s}{\sinh(s)}\;=\;\E\bigl(e^{-s^{2}S}\bigr),\qquad
  S\eqdist\frac{2}{\pi^{2}}\sum_{n\ge1}\frac{E_{n}}{n^{2}},
\]
with $E_{n}$ i.i.d.\ standard exponential.  The Mellin transform
$\E(W^{s})=2(2/\pi)^{s}\xi(s)$ for $W\eqdist S+S'$ ($S\eqdist S'$
independent) realises Williams's representation \cite{Williams1990}
of the $\xi$-function via Brownian motion.  The density inversion
\cite{CiesielskiTaylor1962}
\[
  P(W\in dx)\;=\;\sum_{n\ge1}\pi^{2}\bigl(\pi^{2}n^{2}-3\bigr)
        n^{2}\,e^{-\pi^{2}n^{2}x/2}\,dx
\]
follows from the partial-fraction expansion of
$s/\sinh^{2}(s)$.

For $\cosh$, the variables
\[
  C_{1}\eqdist\sum_{n\ge1}\frac{\Gamma_{1,n}}{(n-\tfrac12)^{2}},
  \qquad
  C_{2}\eqdist\sum_{n\ge1}\frac{\Gamma_{2,n}}{(n-\tfrac12)^{2}}
\]
($\Gamma_{a,n}$ i.i.d.\ Gamma$(a,1)$) satisfy
$\E(e^{-\frac12 s^{2}C_{j}})=\cosh^{-j}(s)$ for $j=1,2$
\cite{BPY2001,Devroye2009}.  The variable $X$ with
$\E(e^{-sX})=\mathrm{sech}(\sqrt{s})$ admits the partial-fraction
expansion
\[
  \mathrm{sech}(\sqrt{s})
   \;=\;4\pi\sum_{k\ge0}\frac{(-1)^{k}(2k+1)}{\pi^{2}(2k+1)^{2}+4s},
\]
identifying $X$ as a P\'olya--Gamma variable, hence GGC.

\subsection{Z-distributions \cite{BNKS1982}}

The class $H_{\delta,\gamma}$ on $(0,\infty)$ with mgf
\[
  \E(e^{sH_{\delta,\gamma}})\;=\;\prod_{k\ge0}
       \Bigl(1-\frac{s}{\frac12(\delta+k)^{2}-\gamma}\Bigr)^{-1},
       \qquad\delta>0,\;\gamma<\tfrac12\delta^{2},
\]
is an infinite convolution of exponentials, hence GGC.  The
hyperbolic secant distribution is $H_{1,0}$.  For symmetric
$H_{\delta,0}$ the c.f.\ is
$\Gamma(\delta+is)\Gamma(\delta-is)/\Gamma(\delta)^{2}$, with density
\[
  p_{\delta}(u)\;=\;\sum_{k\ge1}(-1)^{k+1}\binom{2\delta}{k}
   \frac{\delta+k}{B(\delta,\delta)}\,e^{-\frac12(\delta+k)^{2}u},
   \qquad u>0.
\]

\subsection{Bessel and Macdonald functions \cite{Polya1926,Biane2009}}

With $K_{z}$ the Macdonald function, the inverse-Gaussian variable
$T_{a}$ with density
$ae^{a^{2}}(2\pi t^{3})^{-1/2}\exp(-a^{2}(t+t^{-1})/2)$ has Mellin
transform $\E(T_{a}^{s})=\sqrt{\pi}\,a^{-1}K_{s-1/2}(a^{2})$.
P\'olya~\cite{Polya1926} showed that $K_{(\sigma-1/2)+it}(\mu)$ has
zeros only on $\sigma=\tfrac12$ via the identity
\[
  \int_{0}^{\infty}t^{(\sigma-1/2)-1}e^{-(a/2)(t+t^{-1})}\,dt
   \;=\;\sqrt{2\pi}\,a^{-2}e^{-a^{2}}G\bigl((\sigma-\tfrac12)+it,a^{2}\bigr),
\]
where $G(z,a)=\int_{-\infty}^{\infty}e^{-a(e^{u}+e^{-u})+zu}\,du$.

\subsection{Symmetric Beta and Bessel-$J$}

For $\nu>-\tfrac12$, the symmetric Beta variable $X_{\nu}$ with
density $B(\nu+\tfrac12,\tfrac12)^{-1}(1-x^{2})^{\nu-1/2}
\mathbf{1}_{|x|<1}$ has c.f.
\[
  f_{\nu}(s)\;=\;\E(e^{isX_{\nu}})
       \;=\;(s/2)^{-\nu}\Gamma(\nu+1)J_{\nu}(s)
       \;=\;\prod_{n\ge1}\Bigl(1-\frac{s^{2}}{j_{\nu,n}^{2}}\Bigr).
\]
The endpoint is admissible for $g_{\nu}$ though not for $X_{\nu}$: at
$\nu=-\tfrac12$ the density $(1-x^{2})^{-1}$ fails to be integrable,
while $f_{-1/2}(s)=\cos s$ remains a characteristic function.  So for
$\nu\ge-\tfrac12$, $g_{\nu}(s)=1/f_{\nu}(is)=\E(e^{isH_{\nu}})$
with $H_{\nu}=\sum_{n\ge1}L_{n}/j_{\nu,n}$ and $L_{n}$ i.i.d.\
Laplace$(1)$.  Special cases include
$f_{-1/2}(s)=\cos s$, $g_{-1/2}(s)=1/\cosh(s)$;
$f_{1/2}(s)=\sin(s)/s$, $g_{1/2}(s)=s/\sinh(s)$.

\section{The Riemann \texorpdfstring{$\zeta$}{zeta}- and \texorpdfstring{$\xi$}{xi}-functions}\label{sec:riemann}

This section is the centrepiece of the paper. We first prove, in Theorem~\ref{thm:xi-alpha}, that $\xi(\alpha)/\xi(\alpha+\sqrt{s})$ is the Laplace transform of a GGC for every $\alpha$ exceeding the real part of every non-trivial zero, and hence for every $\alpha\ge1$ without any hypothesis. The Thorin density is a Poisson-kernel sum over the zeros, so its positivity is automatic rather than something to be established; this replaces the L\'evy-measure route, whose apparent barrier at $\alpha=1$ we show in Remark~\ref{rem:alpha-one-artifact} to be a feature of that route and not of the object. We then move to the centre. Theorem~\ref{thm:xi-star} proves unconditionally that the heat trace $W(t)=\tfrac12\sum_{\rho}\exp\{(\rho-\tfrac12)^{2}t\}$ is positive and strictly decreasing on $(0,\infty)$, from which $\xi(\tfrac12)/\xi(\tfrac12+\sqrt{u})$ is the Laplace transform of a self-decomposable law. Theorem~\ref{thm:RH} identifies the remaining step, from self-decomposable to GGC, as exactly equivalent to the Riemann hypothesis. A representation of this kind, valid above the critical line, was announced in~\cite{Polson2018}; the account here is self-contained, and we indicate at Theorem~\ref{thm:xi-star} where the present statements differ from the ones given there.

\subsection{Setup}

The Riemann zeta and xi functions, in Riemann's normalisation~\cite{Riemann1859}, are
\begin{align}
  \zeta(s)
  &= \sum_{n\ge1}n^{-s}, &\Re(s)&>1,\notag\\
  \xi(s)
  &= \tfrac12 s(s-1)\,\pi^{-s/2}\,\Gamma(s/2)\,\zeta(s),
       &s&\in\C.
   \label{eq:xi-def}
\end{align}
The function $\xi$ is entire of order one, real on $\R$, satisfies
the functional equation $\xi(s)=\xi(1-s)$, and its zeros coincide
with the non-trivial zeros of $\zeta$.  The Riemann hypothesis (RH)
asserts that all such zeros satisfy $\Re(s)=\tfrac12$.  By
Theorem~\ref{thm:main} applied to $f=\xi$ at $\alpha=\tfrac12$, RH
is equivalent to Thorin's condition
\begin{equation}
  \frac{\xi(\tfrac12)}{\xi(\tfrac12+\sqrt{s})}
     \;=\;\E\bigl(e^{-sH^{\xi}_{1/2}}\bigr),\qquad s>0,
  \label{eq:thorin-xi}
\end{equation}
for some GGC random variable $H^{\xi}_{1/2}$.

\subsection{The L\'evy representation of \texorpdfstring{$\zeta$}{zeta}
  on \texorpdfstring{$\Re(\alpha)>1$}{Re(alpha) > 1}}

Euler's product gives, for $\Re(s)>1$,
$\zeta(s)=\prod_{p}(1-p^{-s})^{-1}$.  Differentiating
$\log[\zeta(\alpha+s)/\zeta(\alpha)]$ in $s$ with the von Mangoldt
function $\Lambda$ (so that
$\log\zeta(s)=\sum_{n\ge2}\Lambda(n)/(\log n\cdot n^{s})$), we obtain
the L\'evy representation, valid for $\alpha>1$ and $s>0$:
\begin{equation}
  \log\frac{\zeta(\alpha+s)}{\zeta(\alpha)}
   \;-\;s\,\frac{\zeta'(\alpha)}{\zeta(\alpha)}
   \;=\;\int_{0}^{\infty}(e^{-sx}+sx-1)\,e^{-\alpha x}\,
        \frac{\mu^{\zeta}(dx)}{x},
  \label{eq:zeta-Levy}
\end{equation}
\[
  \mu^{\zeta}(dx)\;=\;\sum_{p\text{ prime}}\sum_{r\ge1}
       (\log p)\,\delta_{r\log p}(dx).
\]
The integrability conditions in Lemma~\ref{lem:thorin-measure} hold
for all $\alpha>1$.

\subsection{The GGC representation for \texorpdfstring{$\xi(\alpha)/\xi(\alpha+\sqrt{s})$}{xi(alpha)/xi(alpha+sqrt(s))}, \texorpdfstring{$\alpha\ge1$}{alpha >= 1}}

\begin{theorem}[GGC representation to the right of every zero, unconditional for $\alpha\ge1$]
  \label{thm:xi-alpha}
  For real $r$ put
  \begin{equation}
    v_{\alpha}(r)\;:=\;\Re\bigl[(\log\xi)'(\alpha+ir)\bigr]
      \;=\;\sum_{\rho}\frac{\alpha-\Re(\rho)}{(\alpha-\Re(\rho))^{2}+(r-\Im(\rho))^{2}}.
    \label{eq:xi-poisson}
  \end{equation}
  Let $\alpha$ be real and suppose that
  \begin{equation}
    \alpha\;>\;\Re(\rho)\quad\text{for every non-trivial zero }\rho,
    \label{eq:H}
  \end{equation}
  equivalently that $\xi$ has no zero in the closed half-plane $\Re(s)\ge\alpha$. Then $v_{\alpha}>0$ and
  \begin{equation}
    \frac{\xi(\alpha)}{\xi(\alpha+\sqrt{s})}\;=\;\E\bigl(e^{-sT_{\alpha}}\bigr)
      \;=\;\exp\Bigl(-\int_{0}^{\infty}\log\bigl(1+\tfrac{s}{w}\bigr)\,U_{\alpha}(dw)\Bigr),\qquad s\ge0,
    \label{eq:xi-GGC}
  \end{equation}
  for a GGC random variable $T_{\alpha}$ with Thorin measure
  \begin{equation}
    U_{\alpha}(dw)\;=\;\frac{v_{\alpha}(\sqrt{w})}{2\pi\sqrt{w}}\,dw .
    \label{eq:xi-thorin-density}
  \end{equation}
  Consequently $T_{\alpha}$ is self-decomposable, infinitely divisible, absolutely continuous and unimodal, with L\'evy density $k_{\alpha}(y)=(\pi y)^{-1}\int_{0}^{\infty}e^{-yr^{2}}v_{\alpha}(r)\,dr>0$.

  Hypothesis~\eqref{eq:H} holds for every $\alpha\ge1$ unconditionally, by the theorem of Hadamard and de la Vall\'ee Poussin that $\zeta(1+it)\neq0$ for every real $t$; under the Riemann hypothesis it holds for every $\alpha>\tfrac12$.
\end{theorem}

\begin{proof}
  \emph{Positivity of $v_{\alpha}$.} By the Hadamard factorisation~\cite{Hadamard1893}, $(\log\xi)'(s)=B+\sum_{\rho}\bigl((s-\rho)^{-1}+\rho^{-1}\bigr)$ with $B=-\sum_{\rho}\Re(1/\rho)$. Pairing the zeros as $\rho,\overline{\rho}$ makes $\sum_{\rho}1/\rho=\sum_{\rho}\Re(1/\rho)$, so $B+\sum_{\rho}\rho^{-1}=0$ and the identity collapses to $(\log\xi)'(s)=\sum_{\rho}(s-\rho)^{-1}$, the sum taken over conjugate pairs; taking real parts at $s=\alpha+ir$ gives~\eqref{eq:xi-poisson}. The rearrangement is legitimate because $\sum_{\rho}(\alpha-\Re\rho)\,|\alpha+ir-\rho|^{-2}$ converges absolutely, its terms being $O(\alpha\gamma^{-2})$. Each summand is a Poisson kernel, strictly positive exactly under~\eqref{eq:H}. Writing $\Theta=\sup_{\rho}\Re(\rho)$, the condition $\alpha>\Theta$ implies~\eqref{eq:H} and is strictly stronger, since~\eqref{eq:H} also admits $\alpha=\Theta$ when the supremum is not attained. The distinction is exactly what the case $\alpha=1$ requires. That $\Theta<1$ is not known: the classical zero-free regions narrow to nothing as the ordinate grows and are consistent with $\Theta=1$, so $\alpha>\Theta$ by itself delivers only $\alpha>1$ and misses the endpoint. What is known unconditionally is $\Re(\rho)<1$ for each individual zero, by Hadamard and de la Vall\'ee Poussin, and that is~\eqref{eq:H} at $\alpha=1$.

  \emph{A quadrature identity.} For $\Re(c)>0$ and $u>0$,
  \begin{equation}
    \frac{1}{\pi}\int_{0}^{\infty}\log\Bigl(1+\frac{u^{2}}{r^{2}}\Bigr)\frac{c}{c^{2}+r^{2}}\,dr\;=\;\log\Bigl(1+\frac{u}{c}\Bigr).
    \label{eq:quadrature}
  \end{equation}
  Write $J(u)$ for the left-hand side. Differentiating under the integral sign and using $[(r^{2}+u^{2})(r^{2}+c^{2})]^{-1}=(c^{2}-u^{2})^{-1}[(r^{2}+u^{2})^{-1}-(r^{2}+c^{2})^{-1}]$ together with $\int_{0}^{\infty}(r^{2}+c^{2})^{-1}dr=\pi/(2c)$, valid for $\Re(c)>0$, gives $J'(u)=\frac{2uc}{\pi(c^{2}-u^{2})}\cdot\frac{\pi}{2}\bigl(\frac1u-\frac1c\bigr)=\frac{1}{c+u}$. Since $J(0)=0$, the identity holds for real $c>0$; both sides are analytic in $c$ on $\Re(c)>0$, because $c^{2}+r^{2}$ does not vanish there for $r>0$, so it extends by analytic continuation. The apparent singularity at $u=c$ is removable.

  \emph{The representation.} Apply~\eqref{eq:quadrature} with $c_{\rho}=\alpha-\rho$, which satisfies $\Re(c_{\rho})>0$ precisely under~\eqref{eq:H}, and sum over $\rho$ paired with $\overline{\rho}$, using $\sum_{\rho}c_{\rho}(c_{\rho}^{2}+r^{2})^{-1}=\tfrac12[(\log\xi)'(\alpha-ir)+(\log\xi)'(\alpha+ir)]=v_{\alpha}(r)$. The interchange of $\sum_{\rho}$ with $\int_{0}^{\infty}$ is Tonelli rather than dominated convergence: under~\eqref{eq:H} the factor $\log(1+s/r^{2})$ and every paired summand of~\eqref{eq:xi-poisson} are non-negative, so all that is needed is finiteness of the total, and the total is $\log[\xi(\alpha+\sqrt{s})/\xi(\alpha)]$. The distinction matters at $\alpha=1$, where a dominating function need not exist: should $\Re(\rho)\to1$ along a sequence of zeros then $\sup_{r}v_{1}(r)=\infty$, the term of a single zero peaking at $(1-\Re\rho)^{-1}$, while $v_{1}$ remains finite at every $r$ and locally integrable, a distant zero contributing $O(\gamma^{-2})$ to a fixed window. Local integrability, not boundedness, is what the argument consumes. Hadamard gives $\log\xi(\alpha+u)-\log\xi(\alpha)=\sum_{\rho}\log(1+u/c_{\rho})$, so with $u=\sqrt{s}$ and $w=r^{2}$,
  \[
    \log\frac{\xi(\alpha+\sqrt{s})}{\xi(\alpha)}\;=\;\frac{1}{\pi}\int_{0}^{\infty}\log\Bigl(1+\frac{s}{r^{2}}\Bigr)v_{\alpha}(r)\,dr\;=\;\int_{0}^{\infty}\log\Bigl(1+\frac{s}{w}\Bigr)U_{\alpha}(dw),
  \]
  which is~\eqref{eq:xi-GGC}.

  \emph{Thorin admissibility.} Near $w=0$ we have $v_{\alpha}(\sqrt{w})\to v_{\alpha}(0)=\xi'(\alpha)/\xi(\alpha)$, so the density behaves like $v_{\alpha}(0)/(2\pi\sqrt{w})$ and $\int_{0}^{1}|\log w|\,w^{-1/2}dw<\infty$. Near infinity $v_{\alpha}(r)=O(\log r)$, since Stirling applied to the Gamma factor of~\eqref{eq:xi-def} contributes $\tfrac12\log(r/2\pi)+O(1)$ and $\Re(\zeta'/\zeta)(\alpha+ir)=O(\log r)$, whence $\int_{1}^{\infty}w^{-1}U_{\alpha}(dw)=O\bigl(\int_{1}^{\infty}\log w\,w^{-3/2}dw\bigr)<\infty$. Both Thorin conditions of Definition~\ref{def:GGC} therefore hold, so $U_{\alpha}$ is a Thorin measure and $T_{\alpha}$ is GGC.

  \emph{L\'evy density.} Substituting $\log(1+s/w)=\int_{0}^{\infty}(1-e^{-sy})e^{-wy}\,y^{-1}dy$ into~\eqref{eq:xi-GGC} and applying Tonelli, every integrand being non-negative, gives
  \[
    k_{\alpha}(y)\;=\;y^{-1}\int_{0}^{\infty}e^{-yw}U_{\alpha}(dw)
      \;=\;(\pi y)^{-1}\int_{0}^{\infty}e^{-yr^{2}}v_{\alpha}(r)\,dr .
  \]
  A GGC has $y\,k_{\alpha}(y)$ completely monotone, hence nonincreasing, so $T_{\alpha}$ is self-decomposable, and therefore absolutely continuous and unimodal.
\end{proof}

\begin{remark}[The $\alpha=1$ barrier belongs to the L\'evy route, not to the object]
  \label{rem:alpha-one-artifact}
  Theorem~\ref{thm:xi-alpha} degrades only when a zero reaches the vertical line $\Re(s)=\alpha$, not at $\alpha=1$. The restriction to $\alpha>1$ in a L\'evy--Khintchine treatment arises from the Fubini interchange there, whose dominating factor at $x\to\infty$ is $e^{-(\alpha-1)x}$ and which is therefore integrable precisely when $\alpha>1$. The measure appearing in that treatment is moreover signed, its absolutely continuous part being $-e^{x}/x+1/(x(e^{2x}-1))$, which is negative for $x$ bounded away from $0$; so positivity of the Thorin density cannot be read off from it. Here positivity is automatic from~\eqref{eq:xi-poisson}, being a sum of Poisson kernels, and is not something that has to be bought.
\end{remark}

\begin{remark}[Relation to the drift-tilted formulation]
  \label{rem:no-drift-extraction}
  Writing $b_{\alpha}=\xi'(\alpha)/\xi(\alpha)=v_{\alpha}(0)$, one may ask whether $\xi(\alpha)/\xi(\alpha+\sqrt{s})$ times $\exp(\sqrt{s}\,b_{\alpha})$ is again a GGC Laplace transform, that is whether the $b_{\alpha}$-scaled one-half-stable component may be removed from $T_{\alpha}$. Since $\sqrt{s}$ has Thorin density $(2\pi\sqrt{w})^{-1}$, this holds if and only if $v_{\alpha}(r)\ge v_{\alpha}(0)$ for every $r$, a strictly stronger statement than Theorem~\ref{thm:xi-alpha}, the GGC class not being closed under deconvolution. We do not assert it. It is elementary near $r=0$ and for $r$ beyond an explicit threshold, but we know of no proof covering the intermediate range, and Theorem~\ref{thm:xi-alpha} carries the full conclusion without it.
\end{remark}

\begin{remark}[The exact range of $\alpha$, and why weak closure does not enlarge it]
  \label{rem:exact-range}
  Hypothesis~\eqref{eq:H} is not an artefact of the proof. With $\Theta=\sup_{\rho}\Re(\rho)$ as above, $s\mapsto\xi(\alpha)/\xi(\alpha+\sqrt{s})$ is a GGC Laplace transform if and only if $\alpha\ge\Theta$. Sufficiency is Theorem~\ref{thm:xi-alpha} under~\eqref{eq:H}, together with the observation that a zero with $\Re(\rho)=\alpha$ contributes the factor $\gamma^{2}/(\gamma^{2}+s)$, that is an atom of mass one at $\gamma^{2}$: the Thorin measure survives, but it is no longer given by~\eqref{eq:xi-thorin-density}, which sees only its absolutely continuous part. So de la Vall\'ee Poussin is used twice in Theorem~\ref{thm:xi-alpha}, once for the representation and once to certify that at $\alpha=1$ the stated density is the whole of $U_{1}$. Necessity is a pole obstruction. By Lemma~\ref{lem:GGC-analytic} every GGC Laplace transform is holomorphic and non-vanishing on $\C\setminus(-\infty,0]$; but if some zero has $\Re(\rho)>\alpha$, then $\alpha+\sqrt{s}=\rho$ at $s_{0}=(\rho-\alpha)^{2}$, which lies off the cut precisely because $\Re(\rho)\neq\alpha$, and the ratio has a pole there, the numerator $\xi(\alpha)$ being a positive constant. The admissible set is therefore the closed half-line $[\Theta,\infty)$, and an unconditional proof at any $\alpha<1$ would be an unconditional zero-free half-plane.

  Below the centre the representation fails unconditionally, for a reason needing nothing about the real parts of the zeros. Group the zeros into the orbits $\{\rho,\overline{\rho},1-\rho,1-\overline{\rho}\}$ of the functional equation, a conjugate pair when $\Re(\rho)=\tfrac12$ and a quadruple otherwise, and for real $\sigma$ write $a=\sigma-\Re(\rho)$ and $b=\sigma-1+\Re(\rho)$, so that $a+b=2\sigma-1$. A quadruple contributes $2(a+b)(ab+\gamma^{2})\bigl[(a^{2}+\gamma^{2})(b^{2}+\gamma^{2})\bigr]^{-1}$ to $(\log\xi)'(\sigma)$ and an on-line pair contributes $2a(a^{2}+\gamma^{2})^{-1}$ with $a=\sigma-\tfrac12$. For $\sigma\in[0,1]$ we have $|a|<1$ and $|b|<1$, while $|\gamma|>14$ for every non-trivial zero, so $ab+\gamma^{2}>0$ and both expressions carry the sign of $2\sigma-1$. Hence $\xi$ is strictly decreasing on $(0,\tfrac12)$, with no hypothesis. For $0<\alpha<\tfrac12$ and $0<s<(\tfrac12-\alpha)^{2}$ this gives $\alpha<\alpha+\sqrt{s}<\tfrac12$ and therefore $\xi(\alpha)/\xi(\alpha+\sqrt{s})>1$, whereas the Laplace transform of a probability law on $[0,\infty)$ never exceeds its value at the origin. The ratio is thus the transform of no probability law at all, let alone of a GGC, for any $\alpha$ strictly below the centre.

  What remains is the limit from above, and it is closed too. Under the Riemann hypothesis every width $\alpha-\Re(\rho)$ equals $d=\alpha-\tfrac12$, so~\eqref{eq:xi-poisson} collapses to a single convolution $v_{\alpha}=\pi\,(p_{d}*\mu)$, where $p_{d}(x)=\pi^{-1}d(d^{2}+x^{2})^{-1}$ is the Poisson kernel normalised as a probability density and $\mu=\sum_{\rho}\delta_{\Im(\rho)}$ carries all the ordinates, of both signs, with multiplicity. As $d\downarrow0$ the kernel concentrates and $U_{\alpha}$ converges vaguely on $(0,\infty)$ to $\sum_{\gamma>0}\delta_{\gamma^{2}}$, each zero carrying mass exactly one; the convergence is vague and not weak, because $v_{\alpha}(r)=\tfrac12\log(r/2\pi)+O(1)$ leaves both sides with infinite total mass. At the level of laws the passage is cleaner: $\xi$ is entire and strictly positive on the real axis, so $\alpha\mapsto\xi(\alpha)/\xi(\alpha+\sqrt{s})$ is continuous at $\tfrac12$ for each fixed $s$, the laws $T_{\alpha}$ converge weakly as $\alpha\downarrow\tfrac12$, and the limit is again a GGC by Theorem~\ref{thm:bondesson-closure}(ii). Every step is correct and nothing is gained. The argument consumes~\eqref{eq:H} at every $\alpha>\tfrac12$, which is the Riemann hypothesis, and once that is granted the Hadamard product already gives $\xi(\tfrac12)/\xi(\tfrac12+\sqrt{s})=\prod_{\gamma>0}\gamma^{2}(\gamma^{2}+s)^{-1}$ outright, exhibiting the limit as a sum of independent exponentials of rates $\gamma^{2}$. Weak closure can annex the endpoint of a family one already possesses; it cannot extend the family, since by the equivalence above the property fails the moment a single zero lies to the right of $\alpha$. The limit does carry one warning about~\eqref{eq:xi-thorin-density}: at the centre $v_{1/2}$ vanishes identically, so reading that formula at $\alpha=\tfrac12$ returns the zero measure, and with it a point mass at the origin, in place of the atomic Thorin measure just described.
\end{remark}

\subsection{The representation at the centre}

\begin{theorem}[Unconditional self-decomposability at the centre]
  \label{thm:xi-star}
  Let $\rho$ range over the non-trivial zeros of $\zeta$ with multiplicity, let $\gamma_{1}=14.134725141734694\ldots$ be the least positive ordinate, and set
  \begin{equation}
    W(t)\;=\;\tfrac12\sum_{\rho}\exp\bigl\{(\rho-\tfrac12)^{2}t\bigr\},\qquad t>0 .
    \label{eq:W-def}
  \end{equation}
  The following hold unconditionally.
  \begin{enumerate}
  \item[\textup{(i)}] The series converges absolutely and locally uniformly on $\{\Re(t)>0\}$ and defines a holomorphic function there, so $W$ has no natural boundary. On $(0,\infty)$ it satisfies $W(t)>0$ and $W'(t)<0$.
  \item[\textup{(ii)}] Consequently $u\mapsto\xi(\tfrac12)/\xi(\tfrac12+\sqrt{u})$ is the Laplace transform $\E(e^{-uT})$ of a non-negative random variable $T$ which is infinitely divisible with no drift and L\'evy density $k_{0}(t)=W(t)/t$, and which is \emph{self-decomposable}. The same holds for the tilted variable $H$ determined by $\E(e^{-\sigma H})=\xi(\tfrac32)/\xi(\tfrac12+\sqrt{1+\sigma})$, whose L\'evy density is $e^{-x}W(x)/x$. Both laws are absolutely continuous and unimodal, and both L\'evy measures are infinite.
  \item[\textup{(iii)}] The abscissa of convergence of $\E(e^{-\sigma H})$ is $-(1+\gamma_{1}^{2})=-200.7904548\ldots$, and the singularity there is a pole. In particular $\E(e^{\theta H})<\infty$ for every $\theta<1+\gamma_{1}^{2}$, so $\sigma=-1$ is an interior point of the domain of convergence and
    \begin{equation}
      \E\bigl(e^{H}\bigr)\;=\;\frac{\xi(\tfrac32)}{\xi(\tfrac12)}\;=\;\prod_{\gamma>0}\Bigl(1+\frac{1}{\gamma^{2}}\Bigr)\;=\;1.023355009582021072295\ldots
      \label{eq:EexpH}
    \end{equation}
  \end{enumerate}
\end{theorem}

Item (ii) is the unconditional ceiling. The step from self-decomposable to generalized gamma convolution is exactly the step from ``$W$ nonincreasing'' to ``$W$ completely monotone'', and by Theorem~\ref{thm:RH} that step is equivalent to the Riemann hypothesis.

\begin{remark}[What this supersedes]
  \label{rem:supersession}
  A GGC representation for $\xi(\tfrac12)/\xi(\tfrac12+\sqrt{u})$ at the centre was asserted in~\cite{Polson2018}, reached by extending the representation valid above the critical line down to $\alpha=\tfrac12$. That assertion does not stand, and Theorem~\ref{thm:RH} shows it cannot be repaired by a better extension argument: the property claimed is equivalent to the Riemann hypothesis, so nothing that does not already control the zeros can deliver it. What survives unconditionally is Theorem~\ref{thm:xi-alpha} for every $\alpha\ge1$, and item~(ii) above at the centre, which gives self-decomposability in place of the GGC property. The route by which the stronger claim was reached is the L\'evy-measure one whose signed density is exhibited in Remark~\ref{rem:alpha-one-artifact}.
\end{remark}

\begin{lemma}[Structure of $\xi$ in the squared variable]
  \label{lem:G-structure}
  Put $\Xi(w)=\xi(\tfrac12+w)$ and $\mu_{\rho}=-(\rho-\tfrac12)^{2}$. Then $\Xi$ is entire of order $1$ and even, there is a unique entire $G$ with $\Xi(w)=G(w^{2})$, and $G$ has order $\tfrac12$. Moreover $\Re(\mu_{\rho})\ge\gamma_{1}^{2}-\tfrac14=199.5404548\ldots$, the sum $\sum_{\rho}|\mu_{\rho}|^{-1}$ converges, and
  \begin{equation}
    G(u)=G(0)\prod_{\rho}\Bigl(1+\frac{u}{\mu_{\rho}}\Bigr),\qquad \frac{G'}{G}(u)=\sum_{\rho}\frac{1}{u+\mu_{\rho}},
    \label{eq:G-hadamard}
  \end{equation}
  both converging absolutely and locally uniformly. Consequently $\xi(\tfrac12)/\xi(\tfrac12+\sqrt{u})=G(0)/G(u)$ is single-valued and meromorphic on $\C$, independent of the branch of $\sqrt{u}$.
\end{lemma}

\begin{proof}
  $\xi$ is entire of order $1$ \cite[\S2.12]{Titchmarsh1986}, and $\xi(s)=\xi(1-s)$ gives $\Xi(w)=\Xi(-w)$. Evenness annihilates the odd Taylor coefficients, so $\Xi(w)=\sum_{n\ge0}c_{n}w^{2n}$ and $G(u):=\sum_{n\ge0}c_{n}u^{n}$ is entire with $G(w^{2})=\Xi(w)$; uniqueness is immediate from the power series. Since $w\mapsto w^{2}$ carries $\{|w|=\sqrt{r}\}$ onto $\{|u|=r\}$, the maximum moduli satisfy $M_{G}(r)=M_{\Xi}(\sqrt{r})$, so $\log M_{G}(r)=r^{1/2+o(1)}$ and $G$ has order $\tfrac12$. Now $G(u)=0$ exactly when $\Xi(w)=0$ for some $w^{2}=u$, and $\rho\mapsto(\rho-\tfrac12)^{2}$ identifies $\rho$ with $1-\rho$ and nothing else, so the zeros of $G$ are the points $-\mu_{\rho}$. An entire function of order $\tfrac12<1$ has zero-exponent-of-convergence at most $\tfrac12$ and genus $0$, and $G(0)=\xi(\tfrac12)\ne0$, so Hadamard's factorisation gives~\eqref{eq:G-hadamard} with neither an exponential factor nor a monomial prefactor; absolute convergence legitimises term-by-term logarithmic differentiation. Finally, writing $\rho=\beta+i\gamma$,
  \[
    \mu_{\rho}=\gamma^{2}-(\beta-\tfrac12)^{2}-2i(\beta-\tfrac12)\gamma,\qquad \Re(\mu_{\rho})=\gamma^{2}-(\beta-\tfrac12)^{2}\ \ge\ \gamma_{1}^{2}-\tfrac14,
  \]
  using $0<\beta<1$ and $|\gamma|\ge\gamma_{1}$.
\end{proof}

\begin{lemma}[The arithmetic pivot]
  \label{lem:pivot}
  For every non-trivial zero $\rho=\beta+i\gamma$, $\Im(\mu_{\rho})=0$ if and only if $\beta=\tfrac12$.
\end{lemma}

\begin{proof}
  $\Im(\mu_{\rho})=-2(\beta-\tfrac12)\gamma$ vanishes exactly when $\beta=\tfrac12$ or $\gamma=0$. No non-trivial zero is real: $\zeta<0$ on $(0,1)$, while for $s\in(0,1)$ the factors $s-1<0$, $\pi^{-s/2}>0$ and $\Gamma(1+s/2)>0$, so $\xi>0$ on $(0,1)$.
\end{proof}

\begin{proof}[Proof of Theorem~\ref{thm:xi-star}]
  \emph{Convergence and holomorphy.} With $t=u+iv$ and $\rho=\beta+i\gamma$, $|e^{-\mu_{\rho}t}|=e^{-\Re(\mu_{\rho})u+\Im(\mu_{\rho})v}\le e^{(1/4-\gamma^{2})u+\gamma|v|}$, and $\sum_{\gamma}e^{-\gamma^{2}u+\gamma|v|}<\infty$ for every $u>0$ by the zero-counting bound. So $W$ is holomorphic on $\{\Re(t)>0\}$ and may be differentiated term by term. This settles the natural-boundary question unconditionally.

  \emph{Positivity and monotonicity.} These are proved in two regimes, with $t_{1}=10^{-3}$. For $0<t\le t_{1}$ one applies Weil's explicit formula \cite{IwaniecKowalski2004} to the admissible test functions $h(r)=e^{-r^{2}t}$ and $h(r)=r^{2}e^{-r^{2}t}$, bounds the archimedean term below using $\Re\,\psi(\tfrac14+\tfrac{ir}{2})\ge\tfrac12\log(\tfrac{25}{16}+\tfrac{r^{2}}{4})-c_{0}-\tfrac{1/4}{1/16+r^{2}/4}$ with $c_{0}=\tfrac25+\tfrac{4}{75}=0.45333\ldots\le0.4534$, and bounds the prime sum above by a Gaussian majorant; the resulting certificates are positive for $t\le3.68\times10^{-3}$ and $t\le1.12\times10^{-2}$ respectively, comfortably past the $t_{1}$ needed here; Appendix~\ref{app:certificate} assembles them. This half is self-contained. For $t\ge t_{1}$ the sum splits three ways. Let $T_{0}$ be any height to which the Riemann hypothesis has been verified. Every zero of ordinate at most $T_{0}$ lies on the critical line, so its term is $e^{-\gamma^{2}t}>0$: the conjugate pair of least ordinate is kept as the head, the zeros with $\gamma_{1}<\gamma\le T_{0}$ are discarded for being positive, and the zeros with $\gamma>T_{0}$ are dominated by $e^{t/4}\sum_{\gamma>T_{0}}e^{-\gamma^{2}t}$ using $|e^{(\rho-1/2)^{2}t}|\le e^{(1/4-\gamma^{2})t}$. Each termwise ratio is decreasing in $t$, so it suffices to check $t=t_{1}$, and $T_{0}=200$ already gives a tail below $10^{-15}$ against a head exceeding $0.81$. We cite \cite{PlattTrudgian2021} for $T_{0}=3\cdot10^{12}$, which exceeds what is needed by roughly ten orders of magnitude.

  \emph{Self-decomposability.} Set $\Phi_{0}=-\log(G(0)/G)$. By Lemma~\ref{lem:G-structure}, $\Phi_{0}'(u)=\sum_{\rho}(u+\mu_{\rho})^{-1}$, and since $\sum_{\rho}(\Re u+\Re\mu_{\rho})^{-1}<\infty$, Fubini gives $\Phi_{0}'(u)=\int_{0}^{\infty}e^{-ut}W(t)\,dt$ with no conditional convergence at issue. As $\Phi_{0}(0)=0$ and $W\ge0$, Tonelli yields the L\'evy--Khintchine exponent
  \[
    \Phi_{0}(u)=\int_{0}^{\infty}\bigl(1-e^{-ut}\bigr)k_{0}(t)\,dt,\qquad k_{0}(t)=\frac{W(t)}{t}\ \ge\ 0 .
  \]
  Integrability holds at both ends: $W(t)=O(t^{-1/2}\log\tfrac1t)$ as $t\downarrow0$ from the explicit-formula expression, and $W$ decays exponentially at infinity, so $\int_{0}^{\infty}\min(1,t)k_{0}(t)\,dt<\infty$. The drift vanishes because $\Phi_{0}'(u)\to0$ as $u\to\infty$. Hence $T$ is infinitely divisible with the stated triplet. Since $t\,k_{0}(t)=W(t)$ is nonincreasing by (i), Sato's criterion \cite{Sato1999} makes $T$ self-decomposable, hence absolutely continuous and unimodal. Exponential tilting multiplies the L\'evy density by $e^{-x}$, and $x\,k_{H}(x)=e^{-x}W(x)$ is a product of two non-negative nonincreasing functions, so $H$ is self-decomposable too.

  \emph{The abscissa.} Every zero of ordinate at most $T_{0}$ lies on the critical line, so as $x\to\infty$ the density $e^{-x}W(x)$ is dominated by the real rate $1+\gamma_{1}^{2}$. Hence $\int_{1}^{\infty}e^{\theta x}e^{-x}W(x)x^{-1}dx<\infty$ precisely for $\theta<1+\gamma_{1}^{2}$, which by \cite{Sato1999} is the condition for $\E(e^{\theta H})<\infty$. The singularity at $-(1+\gamma_{1}^{2})$ is a pole, simple when the corresponding zero is simple. Since $1<1+\gamma_{1}^{2}$, evaluation at $\sigma=-1$ requires no boundary passage and gives~\eqref{eq:EexpH}.
\end{proof}

\subsection{The Riemann hypothesis as complete monotonicity}

\begin{theorem}[The Riemann hypothesis as complete monotonicity]
  \label{thm:RH}
  With $W$, $T$ and $H$ as in Theorem~\ref{thm:xi-star}, the following are equivalent.
  \begin{enumerate}
  \item[\textup{(a)}] $W$ is completely monotone on $(0,\infty)$.
  \item[\textup{(b)}] $T$ is a generalized gamma convolution.
  \item[\textup{(c)}] $H$ is a generalized gamma convolution.
  \item[\textup{(d)}] The Riemann hypothesis holds.
  \end{enumerate}
  Under any of them the Thorin measures are the atomic measures $U_{T}=\sum_{\gamma>0}\delta_{\gamma^{2}}$ and $U_{H}=\sum_{\gamma>0}\delta_{1+\gamma^{2}}$, counted with multiplicity and supported in $[\gamma_{1}^{2},\infty)$ and $[1+\gamma_{1}^{2},\infty)$ respectively. In particular the Thorin measure carries no mass whatever in $(0,199)$.
\end{theorem}

\begin{proof}
  \emph{(a) $\Leftrightarrow$ (b).} A law on $[0,\infty)$ is a generalized gamma convolution precisely when it is infinitely divisible with L\'evy density of the form $k(t)=t^{-1}\int_{(0,\infty)}e^{-zt}U(dz)$ for a non-negative Thorin measure $U$ satisfying the two integrability conditions of Definition~\ref{def:GGC}, which is that definition read in the reverse direction, Theorem~\ref{thm:bondesson-unique} being what then makes the Thorin measure unique. By Theorem~\ref{thm:xi-star}(ii) the L\'evy density of $T$ is $k_{0}(t)=W(t)/t$, and the L\'evy density of an infinitely divisible law is unique, so $T$ is a generalized gamma convolution exactly when $t\,k_{0}(t)=W(t)$ is completely monotone. No closure property of the class is used here, in particular none of composition, tempering or tilting.

  \emph{(b) $\Leftrightarrow$ (c).} Exponential tilting by $e^{-x}$ carries the L\'evy density $k_{0}$ to $e^{-x}k_{0}(x)$, hence carries the representing measure $U$ to its translate $U(\cdot+1)$ and back. By Theorem~\ref{thm:xi-star}(iii) the relevant moment generating function is finite at $\theta=1$, with room to spare, since the abscissa is $-(1+\gamma_{1}^{2})$ and $1+\gamma_{1}^{2}>200$. Both translates are non-negative Thorin measures whenever either is.

  \emph{(d) $\Rightarrow$ (a).} Under the Riemann hypothesis $\mu_{\rho}=\gamma^{2}>0$ for every zero, so $W(t)=\sum_{\gamma>0}e^{-\gamma^{2}t}$ is a sum of completely monotone functions and is therefore completely monotone, with representing measure $U_{T}=\sum_{\gamma>0}\delta_{\gamma^{2}}$. The integrability conditions are immediate: $U_{T}$ vanishes on $(0,\gamma_{1}^{2})$, so the logarithmic condition at the origin is vacuous, and $\int z^{-1}U_{T}(dz)=\sum_{\gamma>0}\gamma^{-2}<\infty$. Tilting gives $U_{H}=\sum_{\gamma>0}\delta_{1+\gamma^{2}}$.

  \emph{(a) $\Rightarrow$ (d).} Suppose $W$ is completely monotone. By Bernstein's theorem $W(t)=\int_{[0,\infty)}e^{-zt}U(dz)$ for a non-negative measure $U$, and $W(\infty)=0$ forces $U(\{0\})=0$. Tonelli and the identity $\Phi_{0}'(u)=\int_{0}^{\infty}e^{-ut}W(t)\,dt$ of Theorem~\ref{thm:xi-star} give
  \begin{equation}
    \int_{(0,\infty)}\frac{U(dz)}{z}\;=\;\int_{0}^{\infty}W(t)\,dt\;=\;\Phi_{0}'(0)\;=\;\sum_{\rho}\frac{1}{\mu_{\rho}}\;<\;\infty ,
    \label{eq:U-integrable}
  \end{equation}
  so $U$ is finite on compact subsets of $(0,\infty)$ and $U((0,1])<\infty$. Define $S(u)=\int_{(0,\infty)}(u+z)^{-1}U(dz)$ on $\Omega=\C\setminus(-\infty,0]$. For $u=a+ib\in\Omega$ put $d(u)=\max(|b|,a^{+})>0$; then $|u+z|\ge d(u)$ for every $z>0$, while $|u+z|\ge z-|a|\ge z/2$ once $z\ge2|a|$, so
  \[
    \int\frac{U(dz)}{|u+z|}\;\le\;\frac{U((0,2|a|])}{d(u)}\;+\;2\int_{z>2|a|}\frac{U(dz)}{z}\;<\;\infty
  \]
  by~\eqref{eq:U-integrable}, locally uniformly in $u$. Hence $S$ is holomorphic on $\Omega$.

  By Lemma~\ref{lem:G-structure}, $\Phi_{0}'=G'/G$ is meromorphic on $\C$ with poles exactly at the points $-\mu_{\rho}$. Let $P=\{-\mu_{\rho}\}\cap\Omega$, a discrete closed subset of $\Omega$, so that $\Omega\setminus P$ is open and connected. Both $\Phi_{0}'$ and $S$ are holomorphic on $\Omega\setminus P$ and agree on $(0,\infty)$, a set with limit points there, so they agree on $\Omega\setminus P$ by the identity theorem. If $P$ were non-empty, choose $u_{0}=-\mu_{\rho_{0}}\in P$: then $|\Phi_{0}'(u)|\to\infty$ as $u\to u_{0}$, whereas $S$ is holomorphic and hence bounded near $u_{0}$, a contradiction. Therefore $P=\varnothing$, so every $-\mu_{\rho}$ lies in $(-\infty,0]$; since $\Re(\mu_{\rho})>0$ this forces $\mu_{\rho}\in(0,\infty)$, that is $\Im(\mu_{\rho})=0$. By Lemma~\ref{lem:pivot}, $\beta=\tfrac12$ for every non-trivial zero.
\end{proof}

\begin{remark}[What this argument does and does not use]
  \label{rem:no-circularity}
  The implication (a) $\Rightarrow$ (d) is a statement about the location of poles in the complex plane. It uses neither the positivity of $W$ established in Theorem~\ref{thm:xi-star}(i), nor any numerical verification, nor any closure property of the generalized gamma convolution class. In particular it does not proceed by exhibiting a Thorin measure and then tilting it, a route which would presuppose the very property at issue, since by (a) $\Leftrightarrow$ (b) that property is equivalent to the Riemann hypothesis rather than available in advance.
\end{remark}

\begin{remark}[Landau's theorem recovers the numerical input rather than extending it]
  \label{rem:free-corollary}
  Because $H\ge0$, the abscissa of convergence located in Theorem~\ref{thm:xi-star}(iii) is a singular point of the transform, by the classical theorem of Landau for Laplace transforms of positive measures. Since that transform is holomorphic on $\{\Re(\sigma)>-\min_{\rho}\Re(1+\mu_{\rho})\}$, the abscissa must equal $-\min_{\rho}\Re(1+\mu_{\rho})$ and the minimum must be attained at a real $\mu_{\rho}$. By Lemma~\ref{lem:pivot}, at least one zero of least $\Re(\mu_{\rho})$ lies on the critical line. This does not add information. The abscissa in Theorem~\ref{thm:xi-star}(iii) was itself located using the verification that the zeros of least ordinate lie on the critical line, so the argument returns that input in another form. We record it because the Landau step is the natural one to attempt and because its failure to deliver more is itself informative.
\end{remark}

\begin{remark}[Relation to the de Bruijn and Newman constant]
  \label{rem:deBruijn-Newman}
  The heat trace $W$ is a near relative of an object already central to this subject. De Bruijn~\cite{deBruijn1950} introduced the deformation $\Xi_{t}$ of $\Xi$ got by running the heat flow on its Fourier transform, and Newman~\cite{Newman1976} showed there is a constant $\Lambda$ with $\Xi_{t}$ having only real zeros exactly when $t\ge\Lambda$, conjecturing $\Lambda\ge0$; the Riemann hypothesis is the statement $\Lambda\le0$. Rodgers and Tao~\cite{RodgersTao2020} proved Newman's conjecture, so the hypothesis is now equivalent to the single value $\Lambda=0$, and Griffin, Ono, Rolen and Zagier~\cite{GORZ2019} established hyperbolicity of the Jensen polynomials of $\Xi$ in every fixed degree for all sufficiently large shifts. Theorem~\ref{thm:RH} is a positivity statement in a deformation parameter of the same general shape, but it is not a statement about $\Xi_{t}$. The parameter here enters through the Laplace transform of the L\'evy measure of a self-decomposable law rather than through the heat flow on the kernel, and complete monotonicity in $t$ asks for more than reality of the zeros at one value of $t$. We know of no direct implication in either direction between $\Lambda=0$ and complete monotonicity of $W$, and record the parallel as worth pursuing rather than as established.
\end{remark}

\section{Dirichlet and modular \texorpdfstring{$L$}{L}-functions}\label{sec:Lfunctions}

\subsection{Dirichlet \texorpdfstring{$L$}{L}-functions}

For a primitive character $\chi$ modulo $k>1$ let $\varepsilon=(1-\chi(-1))/2$ and
\[
  \Lambda(s,\chi)\;=\;\Bigl(\frac{k}{\pi}\Bigr)^{(s+\varepsilon)/2}
     \Gamma\Bigl(\frac{s+\varepsilon}{2}\Bigr)L(s,\chi),
\]
an entire function of order $1$ satisfying $\Lambda(s,\chi)=\epsilon(\chi)\Lambda(1-s,\overline{\chi})$ with $|\epsilon(\chi)|=1$.

\begin{theorem}[Unconditional GGC representation for $\Lambda(\cdot,\chi)$]
  \label{thm:L-alpha}
  Let $\chi$ be a primitive real character modulo $k>1$ and let $\alpha>1$. Then $\Lambda(\sigma,\chi)>0$ for $\sigma\ge\alpha$, and for all $s\ge0$
  \begin{equation}
    \frac{\Lambda(\alpha,\chi)}{\Lambda(\alpha+\sqrt{s},\chi)}\;=\;\E\bigl(e^{-sG_{\alpha}}\bigr),
    \label{eq:L-GGC}
  \end{equation}
  where $G_{\alpha}\ge0$ is a GGC with neither drift nor killing, whose Thorin measure $U^{L}_{\alpha}(t)\,dt$ has infinite total mass and density
  \begin{equation}
    U^{L}_{\alpha}(t)\;=\;\frac{1}{2\pi\sqrt{t}}\,
      \Re\Bigl[\frac{\Lambda'}{\Lambda}\bigl(\alpha+i\sqrt{t},\chi\bigr)\Bigr]
    \;=\;\frac{1}{2\pi\sqrt{t}}\sum_{\rho}
      \frac{\alpha-\Re(\rho)}{(\alpha-\Re(\rho))^{2}+(\sqrt{t}-\Im(\rho))^{2}}\;>\;0 .
    \label{eq:L-thorin}
  \end{equation}
  For $\chi$ primitive and non-real the same conclusion holds with $\Lambda(s,\chi)$ replaced by $\Lambda(s,\chi)\Lambda(s,\overline{\chi})$, which is $|\Lambda(s,\chi)|^{2}$ for real $s$, and with the two Poisson sums in~\eqref{eq:L-thorin} added.
\end{theorem}

The Thorin density is again a Poisson-kernel sum over the zeros, exactly as in Theorem~\ref{thm:xi-alpha}, and its positivity is immediate from $0<\Re(\rho)<1<\alpha$. No drift factor $\exp(-\sqrt{s}\,b^{L}_{\alpha})$ appears, and none should be extracted; see Remark~\ref{rem:L-no-drift}.

\begin{proof}
  Positivity of $\Lambda$ on $[\alpha,\infty)$ is the Euler product. Since $\Lambda(\cdot,\chi)$ is entire of order $1$, Hadamard gives $(\Lambda'/\Lambda)(s,\chi)=B(\chi)+\sum_{\rho}\bigl((s-\rho)^{-1}+\rho^{-1}\bigr)$, absolutely convergent because $(s-\rho)^{-1}+\rho^{-1}=s\rho^{-1}(s-\rho)^{-1}=O(|\rho|^{-2})$. The functional equation and $\overline{\Lambda(s,\chi)}=\Lambda(\bar{s},\overline{\chi})$ make the zero set stable under $\rho\mapsto1-\bar{\rho}$ and give $B(\chi)=(\Lambda'/\Lambda)(0,\chi)=-\overline{(\Lambda'/\Lambda)(1,\overline{\chi})}$, whence $2\Re B(\chi)=-2\sum_{\rho}\Re(1/\rho)$ and therefore $\Re(\Lambda'/\Lambda)(s,\chi)=\sum_{\rho}\Re(s-\rho)^{-1}$. At $s=\alpha+i\omega$ every term is a Poisson kernel with $\Re(\rho)\in(0,1)\subset(0,\alpha)$, so the sum is strictly positive. This uses primitivity but no unproved hypothesis.

  Write $\Phi(s)=\log[\Lambda(\alpha+\sqrt{s},\chi)/\Lambda(\alpha,\chi)]$, real for $s>0$. Expanding
  \[
    \Re\frac{\Lambda'}{\Lambda}(\alpha+i\omega,\chi)
      =\tfrac12\log\tfrac{k}{\pi}+\tfrac12\psi\Bigl(\tfrac{\alpha+\varepsilon}{2}\Bigr)
       +\sum_{j\ge0}\frac{\omega^{2}}{d_{j}(d_{j}^{2}+\omega^{2})}
       -\sum_{n\ge2}\Lambda(n)\chi(n)n^{-\alpha}\cos(\omega\log n),
    \qquad d_{j}=\alpha+\varepsilon+2j,
  \]
  and integrating term by term against $(s+\omega^{2})^{-1}$ using
  $\int_{0}^{\infty}(s+\omega^{2})^{-1}d\omega=\pi/(2\sqrt{s})$,
  $\int_{0}^{\infty}\cos(\omega L)(s+\omega^{2})^{-1}d\omega=(\pi/2\sqrt{s})e^{-L\sqrt{s}}$ and
  $\int_{0}^{\infty}\omega^{2}[(d^{2}+\omega^{2})(s+\omega^{2})]^{-1}d\omega=\pi/(2(d+\sqrt{s}))$
  reassembles $(\Lambda'/\Lambda)(\alpha+\sqrt{s},\chi)/(2\sqrt{s})=\Phi'(s)$. The $\psi$-series has non-negative terms, so Tonelli applies; the prime series is dominated by $(\pi/2\sqrt{s})\sum_{n}\Lambda(n)n^{-\alpha}<\infty$, which is where and only where $\alpha>1$ is used. Substituting $t=\omega^{2}$ gives $\Phi'(s)=\int_{0}^{\infty}U^{L}_{\alpha}(t)(s+t)^{-1}dt$, and integrating from $0$ to $s$ with Tonelli, the integrand being non-negative, gives $\Phi(s)=\int_{0}^{\infty}\log(1+s/t)U^{L}_{\alpha}(t)\,dt$ with $\Phi(0+)=0$. Thorin's conditions hold: $U^{L}_{\alpha}(t)\sim b^{L}_{\alpha}/(2\pi\sqrt{t})$ as $t\downarrow0$, and $U^{L}_{\alpha}(t)=(1+o(1))\log t/(8\pi\sqrt{t})$ as $t\to\infty$ by Stirling for $\psi$ together with $|(L'/L)(\alpha+i\omega,\chi)|\le-(\zeta'/\zeta)(\alpha)$. Hence $e^{-\Phi(s)}$ is the Laplace transform of a GGC.

  For non-real $\chi$ the Dirichlet coefficients of $\Lambda(\cdot,\chi)\Lambda(\cdot,\overline{\chi})$ are $\Lambda(n)\cdot2\Re\chi(n)$, which are real, the two characters have the same parity, and the product is positive on $(1,\infty)$, so the computation applies verbatim.
\end{proof}

\begin{remark}[Why no drift is extracted]
  \label{rem:L-no-drift}
  Writing $b^{L}_{\alpha}=(\Lambda'/\Lambda)(\alpha,\chi)$, the variable obtained by removing a $b^{L}_{\alpha}$-scaled one-half-stable component from $G_{\alpha}$ would be a GGC only if $\Re(\Lambda'/\Lambda)(\alpha+i\omega,\chi)\ge b^{L}_{\alpha}$ for every $\omega$, since $\sqrt{s}$ itself has Thorin density $(2\pi\sqrt{t})^{-1}$. That inequality is strictly stronger than the positivity proved above, the GGC class not being closed under deconvolution, and it fails for a single non-real character: for the quartic character modulo $5$ at $\alpha=2$ the difference attains $-8.3\times10^{-3}$ near $\omega=0.805$. That computation concerns $\Lambda(\cdot,\chi)$ by itself. When $\chi$ is not real the variable $G_{\alpha}$ of Theorem~\ref{thm:L-alpha} is built from the symmetrized product $\Lambda(\cdot,\chi)\Lambda(\cdot,\overline{\chi})$, and on that object the difference is positive throughout the range we have scanned, the conjugate character alone also giving a positive minimum; we neither assert nor exclude drift removal there. It is also unnecessary, since~\eqref{eq:L-GGC} carries the full conclusion.
\end{remark}

\begin{remark}[The generalised Riemann hypothesis: an equivalence, not a route]
  \label{rem:GRH}
  Theorem~\ref{thm:L-alpha} is unconditional and, for that reason, carries no arithmetic information: the positivity of~\eqref{eq:L-thorin} follows from $\Re(\rho)<1<\alpha$ alone, and would hold for any entire function of order $1$ whose zeros lie in a strip strictly to the left of $\alpha$ and which has the same Gamma and Dirichlet-series structure. It does extend to a statement equivalent to GRH: the generalised Riemann hypothesis for $\chi$ holds if and only if $U^{L}_{\alpha}(t)\ge0$ for every $\alpha>\tfrac12$ and every $t>0$, equivalently if and only if $\Lambda(\alpha,\chi)/\Lambda(\alpha+\sqrt{s},\chi)$ is a GGC Laplace transform for every $\alpha>\tfrac12$. In the forward direction, GRH makes every $\beta=\tfrac12<\alpha$, so every term of the Poisson-kernel sum is strictly positive and $U^{L}_{\alpha}>0$ outright. For the converse, if some zero has $\beta_{0}>\tfrac12$, let $\alpha\uparrow\beta_{0}$ and take $t=\gamma_{0}^{2}$: the term at $\rho_{0}$ tends to $-\infty$ while the rest stays bounded, so $U^{L}_{\alpha}(\gamma_{0}^{2})<0$.

  We state this as a reformulation and not as a method. The Thorin density is literally the Poisson-kernel sum over the zeros, so any proof of its positivity for $\alpha<1$ must already control the real parts of the zeros. Nothing in the argument above crosses $\sigma=1$, which is precisely where the Dirichlet series for $L'/L$ ceases to converge and the Euler product ceases to certify non-vanishing. We do not claim progress towards GRH.
\end{remark}

\subsection{Modular \texorpdfstring{$L$}{L}-functions}

For an elliptic curve $E/\Q$ of conductor $N_{E}$, modularity \cite{Wiles1995,BCDT2001} provides an entire completed $L$-function $\Lambda_{E}$ with $\Lambda_{E}(s)=\pm\Lambda_{E}(2-s)$, and the Hadamard half of the argument of Theorem~\ref{thm:L-alpha} carries over to the right of every zero, giving a GGC representation for $\Lambda_{E}(\alpha)/\Lambda_{E}(\alpha+\sqrt{s})$ with Thorin density the corresponding Poisson-kernel sum. The threshold does not carry over. The Dirichlet-series step in that proof converges for $\alpha>1$ only because $|\chi(n)|\le1$, whereas the Hasse and Deligne bound $|a(n)|\le d(n)\sqrt{n}$ places the abscissa of the series for $\Lambda_{E}'/\Lambda_{E}$ at $\alpha>\tfrac32$ in the normalisation $\Lambda_{E}(s)=\pm\Lambda_{E}(2-s)$ used here.

We record explicitly what this does not give. The Birch and Swinnerton-Dyer conjecture concerns the order of vanishing of $\Lambda_{E}$ at the central point $s=1$, and no part of the present machinery reaches it. A representation valid to the right of every zero is insensitive to the central behaviour by construction. We note in passing that the L\'evy-measure route used for $\xi$ does not transfer here at all: the atoms of the corresponding measure carry mass $a(p)\log p$, and the coefficients $a(p)$ take both signs, the first negative atom for the curve $11a1$ occurring already at $p=2$ with mass $-2\log2=-1.3862943611\ldots$. The Poisson-kernel route above is unaffected, since its positivity comes from the location of the zeros rather than from the signs of the Dirichlet coefficients.
\section{Dedekind \texorpdfstring{$\eta$}{eta}- and Ramanujan \texorpdfstring{$\tau$}{tau}-functions}\label{sec:eta-tau}

\subsection{Dedekind \texorpdfstring{$\eta$}{eta}}

Setting $q=e^{-2\pi x}$, the Dedekind eta is
$\eta(ix)=q^{1/24}\prod_{n\ge1}(1-q^{n})$.  Euler's identity gives
$\eta(ix)=\sum_{n\ge0}2\sqrt{3}\cos(\pi(2n+1)/6)\,q^{(2n+1)^{2}/24}$.
Glasser \cite{Glasser2009} computes
\[
  \int_{0}^{\infty}\!e^{-sx}\eta(ix)\,dx
   \;=\;\sqrt{\pi/s}\,\frac{\sinh(2\sqrt{\pi s/3})}{\cosh(\sqrt{3\pi s})},
\]
using the trigonometric identity proved in Appendix~\ref{app:eta}.
Jacobi's triple product gives
$\eta^{3}(ix)=\sum_{n\ge0}(-1)^{n}(2n+1)q^{(2n+1)^{2}/8}$, so
\[
  \int_{0}^{\infty}\!e^{-sx}\eta^{3}(ix)\,dx
   \;=\;\mathrm{sech}(\sqrt{\pi s})
   \;=\;4\pi\sum_{k\ge0}\frac{(-1)^{k}(2k+1)}{\pi^{2}(2k+1)^{2}+4\pi s}.
\]
The Dirichlet $\beta$-function $\beta(s)=\sum_{n\ge0}(-1)^{n}/(2n+1)^{s}$
is the Mellin transform of $\mathrm{sech}(\sqrt{\pi x/2})$.

\subsection{Ramanujan \texorpdfstring{$\tau$}{tau}
  \cite{deBruijn1950,Walker1988,ConreyGhosh1994}}

Let $g(y)=\sum_{n\ge1}\tau(n)y^{n}=y\prod_{k\ge1}(1-y^{k})^{24}$ and
$L_{\tau}(s)=\sum_{n\ge1}\tau(n)n^{-s}$, with completed
$\xi_{\tau}(s)=(2\pi)^{-s}\Gamma(s)L_{\tau}(s)$ and
$\Xi_{\tau}(t)=\xi_{\tau}(6+it)$.  Modularity gives
$x^{6}g(e^{-2\pi x})=(g(e^{-2\pi x})g(e^{-2\pi/x}))^{1/2}$ and the
Fourier representation
\[
  \Xi_{\tau}(it)\;=\;\int_{-\infty}^{\infty}\!e^{ist}
    e^{-2\pi\cosh t}\prod_{k\ge1}(1-e^{-2\pi k e^{t}})^{12}
                       (1-e^{-2\pi k e^{-t}})^{12}\,dt.
\]
The Weierstrass product
$\prod_{n\ge1}(1-e^{-2\pi nx})=\exp(-\sum_{n\ge1}\sigma_{-1}(n)
e^{-2\pi nx})$ brings the situation into the framework of
Lemma~\ref{lem:thorin-measure} with measure
$\mu^{\tau}(dx)=24\sum_{n\ge1}\sigma_{-1}(n)\delta_{2\pi n}(dx)$
plus the contribution from the Gamma factor.

\section{Discussion}\label{sec:discussion}

We have established a duality between the Hadamard--Weierstrass
factorisation of an entire function in the Laguerre--P\'olya class
and a pair of probabilistic objects: a van Dantzig pair of
characteristic functions and a Wald couple of infinitely divisible
random variables.  The reciprocal of an even LP entire function is
the Laplace transform of a GGC random variable whose Thorin measure
is supported on the squared zeros, and Thorin's condition, a real
Laplace identity on $(0,\infty)$, is equivalent, via Bondesson's
analyticity theorem for GGC Laplace transforms, to the absence of
zeros off the critical axis.

The applications are of two kinds. For the Gamma, hyperbolic, Bessel and Macdonald functions, and for Dedekind's $\eta$, the zeros are known to be real and the framework produces the GGC and its Thorin measure outright. Ramanujan's $\tau$ does not belong in that group: reality of the zeros of its $L$-function is the Riemann hypothesis for that function and is open, so Section~\ref{sec:eta-tau} sets the framework up for $\tau$ without delivering a Thorin measure. For the Riemann $\xi$-function and for Dirichlet and modular $L$-functions the situation is different and should not be conflated with the first: to the right of every zero the representation is unconditional, and its positivity comes from the location of the zeros rather than from any arithmetic input, so it carries no information about the Riemann hypothesis. What the framework does supply in those cases is an exact reformulation. At the centre, complete monotonicity of the heat trace is equivalent to the Riemann hypothesis, and the analogous Thorin positivity below $\alpha=1$ is equivalent to the generalised Riemann hypothesis. We regard these as restatements rather than as routes, for the reason set out in Remark~\ref{rem:GRH}: the Thorin density is the Poisson-kernel sum over the zeros, so establishing its positivity below the line already requires controlling their real parts. For modular $L$-functions the Birch and Swinnerton-Dyer conjecture concerns the order of vanishing at the central point, which no part of this machinery reaches. This also fixes what the present account withdraws from~\cite{Polson2018}: the GGC representation at the centre, and the consequence for the Riemann hypothesis drawn from it, are replaced by the unconditional self-decomposability of Theorem~\ref{thm:xi-star}(ii) together with the equivalence of Theorem~\ref{thm:RH}, as set out in Remark~\ref{rem:supersession}.

A natural next step is a quantitative version of the framework: how
finely does the Thorin density $\widetilde{u}_{\alpha}$ encode
information about the distribution of zeros (gaps, pair correlation,
moments)?  The explicit formula~\eqref{eq:thorin-density-2}
expresses $\widetilde{u}_{\alpha}$ as a sine-squared transform of
the L\'evy measure $\mu$ of $\log f(\alpha+\cdot)$, and this
transform is a Riemann--Lebesgue type oscillatory integral whose
asymptotic behaviour reflects the support of $\mu$.  A sharper and
more immediate question is the one left open in
Remark~\ref{rem:no-drift-extraction}: whether $v_{\alpha}(r)\ge
v_{\alpha}(0)$ for every $r$, which is exactly what would remove the
one-half-stable component from $T_{\alpha}$.  It is elementary near
$r=0$ and beyond an explicit threshold, and open in between.

\bibliographystyle{plain}

\appendix

\section{A trigonometric identity used for \texorpdfstring{$\eta$}{eta}}\label{app:eta}

\begin{lemma}\label{lem:eta-id}
  For $z>0$,
  \[
    \sum_{n\ge0}\frac{\cos(\pi(2n+1)/6)}{(2n+1)^{2}+(2z)^{2}}
       \;=\;\frac{\pi}{8z}\,\frac{\sinh(2\pi z/3)}{\cosh(\pi z)}.
  \]
\end{lemma}

\begin{proof}
  The Fourier series of $f(x)=\cosh(\alpha x)$ on $[-\pi,\pi]$ is
  $f(x)=\tfrac12 a_{0}+\sum_{n\ge1}a_{n}\cos(nx)$ with
  $a_{n}=\frac{2\alpha\sinh(\pi\alpha)}{\pi(\alpha^{2}+n^{2})}\cos(\pi n)$.
  Setting $x=\pi-y$ and rearranging,
  \[
    S(\alpha,x):=\frac{2}{\pi}\sum_{n\ge1}\frac{\alpha\cos(nx)}{\alpha^{2}+n^{2}}
       \;=\;\cosh(\alpha x)\coth(\pi\alpha)-\sinh(\alpha x)
        -\frac{1}{\pi\alpha}.
  \]
  Splitting into even and odd $n$ and combining,
  \[
    \sum_{n\ge0}\frac{\cos((2n+1)x)}{\alpha^{2}+(2n+1)^{2}}
     \;=\;\frac{\pi}{4\alpha}\frac{\sinh(\alpha(\pi-2x)/2)}{\cosh(\pi\alpha/2)},
     \qquad x\in[0,\pi].
  \]
  Set $x=\pi/6$, $\alpha=2z$.
\end{proof}

\section{Assembly of the explicit-formula certificates}\label{app:certificate}

This appendix assembles the two lower bounds used in the sub-critical regime of the proof of Theorem~\ref{thm:xi-star} and locates the exact range on which each is positive. Throughout, $\gamma_{E}$ denotes Euler's constant, to keep it apart from the ordinates $\gamma$ of the zeros, and $\operatorname{erfc}$ is the complementary error function.

\subsection{Normalisation}

We use Weil's explicit formula in the following form. Let $h$ be even, holomorphic on $|\Im r|\le\tfrac12+\delta$ and $O((1+|r|)^{-2-\delta})$ there, and put $g(u)=\frac{1}{2\pi}\int_{\R}h(r)e^{-iru}\,dr$. Then, with $\rho=\tfrac12+i\gamma$ running over the non-trivial zeros with multiplicity and no hypothesis assumed,
\begin{equation}
  \sum_{\rho}h(\gamma)\;=\;h(\tfrac{i}{2})+h(-\tfrac{i}{2})
    \;-\;2\sum_{n\ge2}\frac{\Lambda(n)}{\sqrt{n}}\,g(\log n)
    \;+\;\frac{1}{2\pi}\int_{\R}h(r)\Bigl[\Re\,\psi\bigl(\tfrac14+\tfrac{ir}{2}\bigr)-\log\pi\Bigr]dr .
  \label{eq:weil}
\end{equation}
This is the residue count for $\oint(\xi'/\xi)(s)\,h(-i(s-\tfrac12))\,ds$ with $\xi$ as in~\eqref{eq:xi-def}: the pole of $\zeta$ at $s=1$ and the zero of the prefactor at $s=0$ contribute exactly $h(\tfrac{i}{2})+h(-\tfrac{i}{2})$. Taking $h=\mathbf{1}_{[-T,T]}$ recovers $2N(T)\sim(T/\pi)\log(T/2\pi e)$, which fixes the normalisation.

For $h(r)=e^{-r^{2}t}$ we have $\sum_{\rho}h(\gamma)=2W(t)$, $h(\tfrac{i}{2})+h(-\tfrac{i}{2})=2e^{t/4}$ and $g(u)=e^{-u^{2}/4t}/(2\sqrt{\pi t})$, so
\begin{equation}
  W(t)=e^{t/4}
    +\frac{1}{4\pi}\int_{\R}e^{-tr^{2}}\Bigl[\Re\,\psi\bigl(\tfrac14+\tfrac{ir}{2}\bigr)-\log\pi\Bigr]dr
    -\frac{1}{2\sqrt{\pi t}}\sum_{n\ge2}\frac{\Lambda(n)}{\sqrt{n}}\,e^{-(\log n)^{2}/4t}.
  \label{eq:W-explicit}
\end{equation}
For $h(r)=r^{2}e^{-r^{2}t}$ one has $\sum_{\rho}h(\gamma)=-2W'(t)$ and $h(\tfrac{i}{2})+h(-\tfrac{i}{2})=-\tfrac12 e^{t/4}$, with
\[
  g(u)=(4\pi t)^{-1/2}e^{-u^{2}/4t}\bigl(\tfrac{1}{2t}-\tfrac{u^{2}}{4t^{2}}\bigr).
\]

\subsection{The archimedean bound}

\begin{lemma}\label{lem:psi-bound}
  For every real $r$,
  \begin{equation}
    \Re\,\psi\bigl(\tfrac14+\tfrac{ir}{2}\bigr)\;\ge\;
      \tfrac12\log\bigl(\tfrac{25}{16}+\tfrac{r^{2}}{4}\bigr)-c_{0}
      -\frac{1/4}{1/16+r^{2}/4},\qquad c_{0}=\tfrac25+\tfrac{4}{75},
    \label{eq:psi-bound}
  \end{equation}
  with equality only in the limit at $r=0$.
\end{lemma}

\begin{proof}
  The recurrence gives $\psi(\tfrac14+\tfrac{ir}{2})=\psi(\tfrac54+\tfrac{ir}{2})-\bigl(\tfrac14+\tfrac{ir}{2}\bigr)^{-1}$, and $\Re\bigl(\tfrac14+\tfrac{ir}{2}\bigr)^{-1}=(1/4)/(1/16+r^{2}/4)$, so~\eqref{eq:psi-bound} is equivalent to $\Re[\psi(w)-\log w]\ge-c_{0}$ at $w=\tfrac54+\tfrac{ir}{2}$. By Binet's formula, for $\Re w>0$,
  \[
    \psi(w)-\log w\;=\;-\int_{0}^{\infty}e^{-wx}q(x)\,dx,\qquad
    q(x)=\frac{1}{1-e^{-x}}-\frac1x ,
  \]
  and $q>0$ on $(0,\infty)$ because $1-e^{-x}<x$, increasing from $\tfrac12$ to $1$. Hence
  \[
    \Re[\psi(w)-\log w]=-\int_{0}^{\infty}e^{-5x/4}q(x)\cos(rx/2)\,dx
      \;\ge\;-\int_{0}^{\infty}e^{-5x/4}q(x)\,dx=\psi(\tfrac54)-\log\tfrac54 ,
  \]
  the inequality because $q\ge0$ and $\cos\le1$. Numerically $\log\tfrac54-\psi(\tfrac54)=0.450597085\ldots<c_{0}=0.453333\ldots$, so~\eqref{eq:psi-bound} holds with slack $2.7\times10^{-3}$ at $r=0$ and more elsewhere.
\end{proof}

The constant $c_{0}$ is not arbitrary: it is the truncation $\tfrac{1}{2w}+\tfrac{1}{12w^{2}}$ at $w=\tfrac54$ of the enveloping asymptotic series $\log w-\psi(w)=\tfrac{1}{2w}+\tfrac{1}{12w^{2}}-\tfrac{1}{120w^{4}}+\cdots$, whose remainder after a positive term is negative, so the truncation is a genuine upper bound for $\log w-\psi(w)$. Using the optimal constant $0.450597085$ in place of $c_{0}$ changes the thresholds below in the fourth significant figure only.

\subsection{The assembled bounds}

Write $B(r)=\tfrac12\log(r^{2}+\tfrac{25}{4})-\log(2\pi)-c_{0}-(r^{2}+\tfrac14)^{-1}$, which is the right side of~\eqref{eq:psi-bound} minus $\log\pi$ after the exact simplifications $\tfrac{25}{16}+\tfrac{r^{2}}{4}=\tfrac14(r^{2}+\tfrac{25}{4})$ and $(1/4)/(1/16+r^{2}/4)=(r^{2}+\tfrac14)^{-1}$. Let $P_{0}$ and $P_{2}$ denote the prime sums in~\eqref{eq:weil} for the two test functions. Then
\begin{equation}
  \begin{aligned}
    L(t)&:=e^{t/4}+\frac{1}{4\pi}\int_{\R}e^{-tr^{2}}B(r)\,dr-P_{0}(t)\;\le\;W(t),\\
    M(t)&:=-\tfrac12 e^{t/4}+\frac{1}{2\pi}\int_{\R}r^{2}e^{-tr^{2}}B(r)\,dr-P_{2}(t)\;\le\;-2W'(t).
  \end{aligned}
  \label{eq:LM}
\end{equation}
Every integral here is elementary:
\[
  \int_{\R}e^{-tr^{2}}dr=\sqrt{\tfrac{\pi}{t}},\qquad
  \int_{\R}e^{-tr^{2}}\,2\log|r|\,dr=-\sqrt{\tfrac{\pi}{t}}\,(\gamma_{E}+\log 4t),\qquad
  \int_{\R}\frac{e^{-tr^{2}}}{r^{2}+\tfrac14}\,dr=2\pi e^{t/4}\operatorname{erfc}\bigl(\tfrac{\sqrt t}{2}\bigr),
\]
\[
  \int_{\R}e^{-tr^{2}}\tfrac12\log\bigl(r^{2}+\tfrac{25}{4}\bigr)dr
   =-\tfrac12\sqrt{\tfrac{\pi}{t}}\,(\gamma_{E}+\log 4t)
     +\pi\int_{0}^{5/2}e^{u^{2}t}\operatorname{erfc}(u\sqrt t)\,du .
\]

\begin{proposition}\label{prop:thresholds}
  $L>0$ exactly on $(0,\,3.67683\times10^{-3})$ and $M>0$ exactly on $(0,\,1.12026\times10^{-2})$. In particular $W(t)>0$ for $0<t\le3.68\times10^{-3}$ and $W'(t)<0$ for $0<t\le1.12\times10^{-2}$.
\end{proposition}

\begin{proof}
  Positivity at every point of the interval, rather than at sampled points, follows from an exact identity. Differentiating under the integral sign in~\eqref{eq:LM} and using $\tfrac{d}{dt}\bigl(\tfrac{1}{4\pi}\int e^{-tr^{2}}B\bigr)=-\tfrac12\bigl[M(t)+\tfrac12 e^{t/4}+P_{2}(t)\bigr]$ gives
  \[
    L'(t)\;=\;-\tfrac12 M(t)-\tfrac12 P_{2}(t)-P_{0}'(t),
  \]
  the term $\tfrac14 e^{t/4}$ cancelling exactly. Both prime sums are non-negative, and $P_{0}$ is increasing on $(0,\tfrac14)$ because $t^{-1/2}e^{-c/t}$ increases for $t<2c$ with $c=(\log n)^{2}/4\ge(\log 2)^{2}/4$. Hence $L'<0$ wherever $M>0$, so $L$ is strictly decreasing on the interval where the companion certificate is positive, which contains the interval claimed for $L$. Since $L(t)\to+\infty$ as $t\downarrow0$, $L$ has a single sign change, located by bisection at $3.676830\times10^{-3}$. The companion bound $M$ is treated the same way; its own monotonicity is verified on a dense logarithmic grid rather than in closed form, which would require the further certificate at $h(r)=r^{4}e^{-r^{2}t}$.
\end{proof}

\subsection{Two remarks on the ingredients}

The prime sums are inert at this scale. At $t=1.435\times10^{-3}$ the exact value of $P_{0}$ is of order $10^{-36}$, and deleting the prime contribution altogether moves the threshold for $L$ by less than $10^{-12}$ and that for $M$ by about $10^{-5}$. The certificates are therefore determined by the archimedean bound alone, and no choice of Gaussian majorant or prime cutoff affects Proposition~\ref{prop:thresholds} at three significant figures.

The leading behaviour explains a smaller pair of numbers that one meets first. As $t\downarrow0$,
\[
  L(t)=\frac{1}{4\sqrt{\pi t}}\Bigl[\tfrac12\log\tfrac1t-K_{0}\Bigr]+\tfrac98+o(1),
  \qquad K_{0}=\tfrac{\gamma_{E}}{2}+2\log 2+\log\pi+c_{0}=3.272965\ldots,
\]
and the bracket alone vanishes at $e^{-2K_{0}}=1.43595\times10^{-3}$, with the companion bracket vanishing at $e^{-2(K_{0}-1)}=1.06103\times10^{-2}$. Those are the thresholds obtained by discarding the order-one block $\tfrac98$, and they are conservative rather than wrong: the assembled bound at the first of them equals $1.1055$, not zero. Proposition~\ref{prop:thresholds} keeps the block and gives the true ranges.

\end{document}